\documentclass{amsart}
\usepackage{amsfonts,amssymb}

\def\Fscr{{\mathcal F}}

\def\Mscr{{\mathcal M}}

\def\QQ{{\mathbb Q}}

\def \PP{{\mathbb P}}

\def \ZZ{{\mathbb Z}}

\def \NN{{\mathbb N}}

\def\spec{\operatorname {Spec}}
\def\SL{\operatorname {SL}}
\def\Gr{\operatorname {Gr}}
\def\maxspec{\operatorname {maxSpec}}

\def\Spec{\operatorname {Spec}}

\def\Ext{\operatorname {Ext}}
\def\Hom{\operatorname {Hom}}

\def\Supp{\operatorname {Supp}}
\def\sup{\operatorname {sup}}

\def\im{\operatorname {im}}

\def\Tor{\operatorname {Tor}\nolimits}
\def\End{\operatorname {End}}

\def\r{\rightarrow}

\def\xs{x_1, x_2, \dots, x_d}

\def\Ext{\operatorname {Ext}\nolimits}
\def\rad{\operatorname {rad}\nolimits}
\def\rk{\operatorname {rk}\nolimits}
\def\Hom{\operatorname {Hom}\nolimits}

\def\im{\operatorname {im}\nolimits}
\def\End{\operatorname {End}\nolimits}

\def\trace{\operatorname{trace}\nolimits}
\numberwithin{equation}{section}
\def\pretendhaswidth#1#2{\setbox0\hbox{$#1$} \hbox to \wd0{\hss $#2$ \hss}}

\newtheorem{lemma}{Lemma}[section]

\newtheorem{theorem}[lemma]{Theorem}

\newtheorem{lemmas}{Lemma}[subsection]
\newtheorem{propositions}[lemmas]{Proposition}
\newtheorem{theorems}[lemmas]{Theorem}
\newtheorem{corollarys}[lemmas]{Corollary}

\theoremstyle{definition}

\newtheorem{conjecture}[lemma]{Conjecture}

\newtheorem{examples}[lemmas]{Example}
\newtheorem{definitions}[lemmas]{Definition}

\newtheorem{questions}[lemmas]{Question}

\newtheorem{step}{Step}

\theoremstyle{remark}

\newtheorem{remarks}[lemmas]{Remark}

\newtheorem{observations}[lemmas]{Observation}

\title{Simplicity of rings of differential operators in
prime characteristic}
\keywords{Rings of differential operators, Strong F-regularity, Finite
representation type}
\subjclass{Primary 16S32, 16G60, 13A35} 
\author{Karen E. Smith}
\address{Math Dept. 2-167 \\
Massachusetts Institute of Technology\\ 02139 Cambridge
\\Massachusetts}
\email{kesmith@math.mit.edu}
\thanks{The first author is supported by the  US National Science Foundation.}
\author{Michel Van den Bergh}
\address{Limburgs Universitair Centrum\\ Departement WNI\\ Universitaire
Campus\\ 3590 Die\-pen\-beek \\ Belgium}
\thanks{The second author is a senior researcher at the
NFWO}
\email{vdbergh@luc.ac.be, 
http://www.luc.ac.be/Research/Algebra/}
\begin{document}
\begin{abstract} Let $W$ be a finite dimensional representation of a linearly 
reductive group $G$ over a field $k$. Motivated by their work on classical 
rings of invariants, Levasseur and Stafford asked whether the
ring of invariants under $G$ of the symmetric algebra of $W$ has a simple ring
of differential operators.
In this paper, we show that this is true in prime characteristic.
Indeed, if $R$ is a graded subring of a polynomial ring over a perfect field 
of characteristic $p>0$ and if the inclusion $R\hookrightarrow S$ splits, then
$D_k(R)$ is a simple ring. In the last section of the paper, we discuss how
one might try to deduce the characteristic zero case from this result.
As yet, however, this is a subtle problem and the answer to the question of
Levasseur and Stafford remains open in characteristic zero.
\end{abstract}
\maketitle
\tableofcontents
\section{Introduction}
\label{sec1}
The ring  
of differential operators on the coordinate ring $R$ of a smooth 
connected affine  algebraic variety over a field $k$ is  fairly 
well understood.
 This ring, denoted $D_k(R)$, is a simple ring; in the event that $k$ has
characteristic zero, $D_k(R)$ is a finitely generated $k$-algebra.
Despite intense study by various authors, however, the
 structure of $D_k(R)$ for a non-smooth variety remains a mystery. 
 It is known that it need  not be simple 
nor finitely generated in general \cite{BGG2}.
A problem remains to identify classes of varieties for which the corresponding
ring of differential operators is simple or is finitely generated.

A general expectation persists that rings of differential operators on rings 
 of invariants for linearly
 reductive groups should have various nice properties \cite{LSt}
 \cite{Musson}\cite{Schwarz}\cite{VdB5}.
One such expectation is reflected in the following conjecture.
\begin{conjecture}\label{con1}
 Let $W$ be a finite dimensional $k$-representation of a linearly reductive
group $G$. Let $R$ be the ring of invariants for $G$ acting on  the 
symmetric algebra $S(W)$. Then $D_k(R)$ is a simple ring.
\end{conjecture}
This conjecture was stated as a question by Levasseur and Stafford 
\cite{LSt}
 in the slightly more restrictive
setting of  $G$ reductive  and $k$ of characteristic zero.
To date, much research into the structure of rings of differential
operators has focused on the 
case of affine rings over a field of characteristic zero.
This may be founded in a (perhaps misguided) feeling that 
rings of differential operators are better behaved in 
 characteristic zero.
In fact, differential operators have an especially 
interesting form in prime characteristic,
and their nice structure enables us to prove a general theorem 
that implies Conjecture \ref{con1}
in characteristic $p$.

In a fairly delicate argument
using the structure theory of primitive ideals in an enveloping algebra,
 Levasseur and Stafford were able to prove Conjecture \ref{con1} for the
``classical'' rings of invariants in characteristic zero \cite{LSt}.  
Earlier, Kantor \cite{Kantor} had
 proven the conjecture for finite groups, assuming the ground field
to be the complex numbers.
 Van den Bergh
has studied some representations of $\SL(2, k)$ when $k$ has characteristic 
zero \cite {VdB13}.  The case of torus invariants in 
characteristic zero is covered by \cite{MVdB}. 
Very little had been known about the structure of $D_k(R)$, however, in the
case where $k$ has characteristic $p>0$.

\medskip
The simplicity of $D_k(R)$ imposes strong restrictions on $R$.
If $D_k(R)$ is a simple ring, then $R$ must be simple as a $D_k(R)$
module (see \S\ref{sec5}). 
For any reduced ring $R$, it is easy to see that
each minimal prime of $R$ is a $D_k(R)$ submodule of $R$, so 
 that simplicity of 
$D_{k}(R)$ implies that $R$ is a domain.
In \cite{VdB5}, it is shown that simplicity of $D_k(R)$ implies 
that $R$ is Cohen-Macaulay.
In characteristic $p>0$, 
the simplicity of $D_k(R)$ forces the tight closure of an ideal 
$I$, usually a very subtle and difficult to compute closure operation,
to take an especially simple form
\cite{smith}. 

\medskip

When $R$ is the ring of invariants of a linearly reductive group 
acting on the symmetric algebra $S = S(W)$ 
 of a finite dimensional $k$-representation $W$,
the Reynolds operator provides an $R$ module splitting
of the inclusion map $R \hookrightarrow S$.
The simplicity of $D_k(R)$  may be the result of
 the 
 elementary algebraic features of $R$ inherited by virtue
of its being a direct summand of a regular ring 
rather than a consequence of 
 some of the  more subtle issues arising because of the structure and action
 of $G$.  We propose the following conjecture.

\begin{conjecture}
Let $R \hookrightarrow S$ be an inclusion of $k$-algebras, where $k$ is
any field. Assume that this inclusion splits as a map of $R$ modules.
If $D_k(S)$ is a simple ring, then $D_k(R)$ is a simple ring.
\end{conjecture}

In the question of   Levasseur and  Stafford,
$S$ is a polynomial ring and $R$ is some graded subring.  
Thus to prove Conjecture \ref{con1}, it would suffice to prove
this more general  conjecture 
even just for the case where $R$ is a graded subring of a polynomial ring $S$.
This is exactly what we accomplish in this paper in prime characteristic.
The following 
 important corollary  follows easily from  the main theorem of this paper 
(Theorem \ref{th41}).
\begin{theorem}
Let $R$ be a graded subring of a polynomial ring $S$ over a 
perfect field $k$ of characteristic $p>0$. Assume that the inclusion 
$R \hookrightarrow S$ splits in the category of graded $R$ modules.
Then $D_k(R)$ is a simple ring.
\end{theorem}

With the notation of the theorem above,
 $D_k(R)$ 
is simply the ring of all additive maps from $R$ to itself linear
over some subring $R^{p^e}$ of ${p^e}$-th powers of the 
elements of $R$ (see \S\ref{diffopsp}). Our work
heavily exploits this very interesting description of the ring of
differential operators.

\bigskip

The study of differential operators in characteristic $p>0$ 
is partially motivated by the connections with the theory of
\emph{tight closure}. 
Tight closure, introduced by 
M. Hochster and C. Huneke,
 is a closure operation performed on ideals in a ring
of prime characteristic which has led to deep new insight into 
the structure of commutative rings containing a field. We refer the reader
to \cite{HH}  for more about tight closure. 
The relationship with differential operators is developed by the
first author in \cite{smith}.
Although this paper is  independent of the
 theory of tight closure, some of  
the motivation for this work,
and many of the ideas within, 
are inspired by tight closure.

\medskip

The paper is organized as follows. Section \ref{sec2} establishes the
notation and summarizes the relevant facts 
to  be used throughout the paper.
Section \ref{sec3}   contains a study of  the behavior  of modules
under the Frobenius functor.  In particular, we
introduce the concept of ``Finite F-representation type.''
This is a fairly strong
 representation theoretic property of commutative rings of characteristic 
$p>0$. Its name is intended to recall the
 similar (though much stronger)  property
of finite representation type. 
In Section \ref{sec4}, the main results about the
simplicity of rings of differential operators in characteristic $p>0$
are developed. 
 In the final section 
of this paper, we address the question of 
how one may attempt to deduce the characteristic zero case from this
theorem. As yet, however, this problem is a subtle one and the above 
conjectures
 remain open in  characteristic zero.

The authors thank Mel Hochster for illuminating various aspects
of \cite{HH3}; his comments regarding the behavior of discriminants 
were particularly helpful in straightening out the proof of 
Theorem \ref{Dsimplered}. The first author also expresses warm thanks to
Ragnar Buchweitz and Frank Schreyer for an engaging discussion 
about Atiyah's classification of vector bundles on an elliptic curve
that was helpful in understanding \cite{Tango}.

\section{Generalities}
\label{sec2}

Throughout this paper, the word ``module'' always means
 left module, except where
otherwise indicated. The notation $M_R$ indicates that $M$ is a right 
$R$-module, and the notation ${}_RM$ indicates that $M$ is a left $R$-module.
Frequently, an abelian group $M$
 may be a module or bi-module over several different rings,
in which case the notation $M_R$  (or ${}_RM$) 
indicates that the  structure under consideration is the right (or left)
$R$-module structure.

\subsection{Differential Operators: Definitions}
\label{diffops}
Let $k$ be any commutative ring and let $R$ be any commutative $k$-algebra.
For any 
two $R$ modules $M$ and $N$, the module
$D_{R/k}(M, N)$ (or $D_{k}(M, N)$ when $R$ is understood)
   of $k$ linear 
 differential operators from $M$ to $N$ is a certain
distinguished submodule  of $\Hom_{k}(M, N)$, defined inductively
as follows:
$$
D_k(M, N) = \bigcup_{n \in \mathbb N}D^n_k(M, N),
$$
where 
$$
D^0_k(M, N) = \Hom_R(M, N) {\text{   and}}
$$
$$
D^n_k(M, N) = \{\theta \in \Hom_k(M,N) | {\text{  for all }} r \in R,
[r, \theta] \in D^{n-1}_k(M, N)\}.
$$
The symbol $[r, \theta]$ denotes the commutator operator
$\,(r \circ \theta - \theta \circ r)\, \in \Hom_k(M, N)$, where
the symbol $``\circ$'' denotes composition of operators.
A differential operator $\theta \in D^n_k(M, N)$ but not in $D_k^{n-1}(M, N)$ 
is said to be of {\it order n.\/} 
When $M = N$, the differential operators $D_k(M, M)$ form a ring,
denoted $D_k(M)$.
Note that $D_k(M, N)$ is a $D_k(N)-D_k(M)$ bimodule where the
action is given by composition of maps.
 Of particular interest is the ring of differential 
operators on $R$ itself, $D_k(R)$.  
The unadorned symbol $D(R)$ indicates that the base ring $k$ is $\ZZ$.

Another, entirely equivalent, point of view
is the following. Consider the ring $R\otimes_k R$.
Note that $\Hom_k (M, N) $ is a (left)  module over 
$R\otimes_k R$ is an obvious way: an element $r \otimes s$ acts on 
$\theta$ to produce $r \circ \theta \circ s \in \Hom_k (M, N) $.
Let $J_{R/k}$ denote the kernel of the multiplication map
 $R \otimes_k R \r R$ sending  an element $r \otimes s$ to the product $rs$.
Note that $J_{R/k}$ is generated by 
elements of the form $r \otimes 1 - 1 \otimes r$,
and that such an element acts on $\theta \in \Hom_k (M, N)$ to produce the
commutator $[r, \theta]$.  
It is thus clear that an operator $\theta \in \Hom_k (M, N)$ is a differential 
operator of order less than or equal to $n$ if  and only if 
$\theta $ is annihilated by $J^{n+1}$.
That is, there are canonical isomorphisms
$$
D^n_k(M, N) = \Hom_{R\otimes_k R}((R\otimes_k R)/J_{R/k}^{n+1}, \Hom_{k}(M, N))
= \Hom_R (P^n_{R/k}\otimes M, N)
$$
where $P^n_{R/k}$ is simply  $(R\otimes_k R)/J^{n+1}_{R/k}$,
regarded as an $R$-bimodule, and the homomorphisms are as left $R$
modules. By definition, an
 element $m$ of an $S$-module $\Mscr$ is in the zero-th local cohomology
 module $H^0_J(\Mscr)$ with support in the ideal $J \subset S$ if
$m$ is annihilated by some power of $J$.
Thus,  differential operators can be thought of as the elements of the
local cohomology module
\begin{equation}
\label{localcom}
D_k(M, N) = H^0_{J_{R/k}}(\Hom_k (M, N)).
\end{equation}
\subsection{Derived functors of differential operators}
\label{ssec215}
Of substantial interest to us are the higher derived functors of differential
operators. This is because, as we will show in Section \ref{sec6}, they 
essentially control the behavior of differential operators under reduction to 
prime characteristic. Denote by $R^iD_k(M,-)$ the right derived functors 
of the left exact covariant functor $D_k(M,-)$. Clearly
\[
R^iD_k(M,N)=\varinjlim_n \Ext^i_R(P^n_{R/k}\otimes_R M,N)
\]
Although $D_k(-, N)$ is a left exact
contravariant functor (so that it would  make
 sense to compute its derived functors), in general
 $R^iD_k(M,N)$ is not the derived functor of its first 
argument; indeed, there is no reason why $R^iD_k(R,M)$ should vanish for $i>0$
(see example \ref{example51}). However, 
using the ``associativity formula'' 
$\Hom_R(M, \Hom_R(P^n_{R/k}, N)) =  \Hom_R(P^n_{R/k}\otimes M, N),$
we do see that there are two spectral sequences with the same limit whose
$E^2$ terms are (dropping the subscripts from the notation):
$$
\Ext^p(M, \Ext^q(P^n, N)) \qquad 
{\text{  and  }}\qquad \Ext^q(\Tor_p(P^{n}, M), N)
$$
(see chapter XVI of \cite{ce}).

In very good cases, these spectral sequences collapse, as we now describe.
Recall that 
 a $k$-algebra $R$ is \emph{formally smooth}
if,  for every $k$-algebra $C$ and every ideal $J \subset C$  such that
$J^2 = 0$, the map $\Hom_k(R, C) \rightarrow \Hom_k(R, C/J)$ is
surjective (\cite{EGA}, IV.0.19.3.1, or IV.17.1.1). 
If $k$ is a field, then the coordinate ring $R$ 
 of  a  smooth affine algebraic variety  over $k$ is $k$-formally smooth, 
as is any local
 ring of a smooth $k$-variety.

\begin{propositions}
\label{CCformsooth}
 Assume that $R$ is formally smooth over $k$, 
 and that each $P^n_{R/k}$ is a 
finitely presented (left) $R$ module. 
Then 
\begin{enumerate}
\item  $R^iD_k(R,N)=0$ for $i>0$;
\item $R^iD_k(M,N)$ is a derived functor of its first argument;
\item  Assume in addition that $R$ is noetherian and that $M$ is
  finitely generated. Then there are natural identifications
\begin{align*}
R^iD_k(M,N)&=\Ext^i_R(M,(N\otimes_R D_k(R))_R)\\
&=\Ext^i_{D_k(R)}(M\otimes_R D_k(R),N\otimes_R D_k(R))
\end{align*}
\end{enumerate}
\end{propositions}
\begin{proof}
The main point is that if 
$R$ is formally  smooth over $k$,  then each  $P^n_{R/k}$ is left and right 
projective as an $R$ module (\cite{EGA}, 16.10.1). It easily follows that
$$
\Ext^i(M, \Hom(P^n, N)) \cong \Ext^i(P^{n} \otimes M, N).
$$  
for all $i$ (this can be seen also from the collapsing of both
spectral sequences at $E^2$).
 Taking the direct limit (which commutes with the
computation of cohomology),  we see that both (1) and (2) hold.

For the first identification in  part (3), 
we have
\begin{align*}
R^iD_k(M, N) = & \varinjlim \Ext^i_R(M, \Hom_R(P^n, N))\\
\cong & \Ext^i_R(M, \varinjlim \Hom_R(P^n, N))\\
\cong & \Ext^i_R(M,   (N \otimes_R D(R))_R).\\
\end{align*}
The first isomorphism above follows from the fact that $M$ is finitely
generated and  $R$ is noetherian. Under these hypotheses it is
easily seen that $\Ext$ commutes with direct limits in it second
argument.

The last isomorphism uses the assumption that each $P^n$ is
finitely presented and projective, in order to verify that the natural
map $N \otimes \Hom_R(P^n, R) \rightarrow \Hom_R(P^n, N)$ is an
isomorphism. Taking the direct limit, we have the natural isomorphism
of right $D(R)$ modules (and hence right $R$ modules) $N \otimes D(R)
= D(R, N)$.

For the second identification in (3), 
we point out that $D(R)$ is a flat $R$ module, since it is a
direct limit
of projective $R$ modules. Furthermore, for any right $D(R)$ module $I$,
and any right $R$ module $M$, we have the natural adjointness
$$
\Hom_R(M, I_R) = \Hom_{D(R)}(M \otimes_R D(R), I).
$$
It follows that
 any injective right $D(R)$ module is also injective considered as
a right $R$ module, since the functor 
$\Hom_R( - , I_R) = \Hom_{D(R)}( - \otimes_R D(R), I)$ is exact.
The computation of  $\Ext^i_R(M, (N\otimes_R D(R))_R)$ can therefore
be accomplished by resolving $N\otimes D(R)$ by right $D(R)$
injectives,  and then viewing this resolution as a resolution by
right $R$ modules.  In light 
 of the adjointness
above, the identification with $\Ext^i_{D_k(R)}(M\otimes_R
D_k(R),N\otimes_R D_k(R))$ is immediate.
\end{proof}

The formula \eqref{localcom} suggests a corresponding formula for higher local
cohomology. The following proposition shows that this is essentially correct.
\begin{propositions}
\label{CCformgeneral}
Assume that $R$ and $M$ are projective over $k$ and that $M$ is 
a finitely
generated $R$-module. Assume furthermore that $R\otimes_k R$ is 
Noetherian\footnote{$R$ is a quotient of $R\otimes_k R$, so if $R\otimes_k R$
is Noetherian, then so is $R$. The converse is false; a counter example
is given by a rational function field in infinitely many 
variables.}. Then 
\[
R^iD_k(M,N)=H^i_{J_{R/k}}(\Hom_k(M,N))
\]
\end{propositions} 
\begin{proof}
Because $M$ is a projective $k$ module, the functor 
$H^0_{J_{R/k}}(\Hom_k(M,-))$ is a left exact  covariant functor, and 
the $\{H^i_{J_{R/k}}(\Hom_k(M,-))\}$ form a $\delta$ functor (so that 
short exact sequences of $R$ modules give rise to 
 long exact sequences involving $H^i_{J_{R/k}}(\Hom_k(M,-))$).
  By \eqref{localcom}, $H^i_{J_{R/k}}(\Hom_k(M,-))$ 
 agrees with $R^iD(M, N)$ when  $i=0$.
 Therefore, to prove the proposition, it will
 be sufficient to show that  $H^i_{J_{R/k}}(\Hom_k(M, -))$ vanishes
 on injective $R$ modules when $i > 0$.
Now
\[
H^i_{J_{R/k}}(\Hom_k(M,I))=\varinjlim_n\Ext^i_{R\otimes_k R}
(P^n,\Hom_k(M,I))
\]
Let $A^\cdot$ be a projective $R\otimes_k R$-resolution of 
$P^n$. 
We have
\begin{align*}
\Ext^i_{R\otimes_k R}
(P^n,\Hom_k(M,I)) &= H^i(\Hom_{R\otimes_k R}
(A^\cdot, \Hom_k(M,I)))\\
&=H^i(\Hom_R(A^\cdot\otimes_R M,I))\\
&=\Hom_R(H_i(A^\cdot\otimes_R M),I)
\end{align*}
The assumption that  $R$ is $k$-projective implies that  any projective
$R\otimes_k R$ module is also a projective $R$ module (on either side).
Thus the $R\otimes R$ projective resolution $A^\cdot$  can be viewed
as 
a projective resolution of
$P^n$ as a right $R$-module.
So
\[
H_i(A^\cdot\otimes_k M)=\Tor_i^R((P^n)_R, M)
\]
and we obtain
\begin{equation}
\label{CCbasicform}
H^i_{J_{R/k}}(\Hom_k(M,I))=\varinjlim_n 
\Hom_R(\Tor^R_i(P^n,M),I)
\end{equation}
Now let $B^\cdot$ be a resolution of $M$ by finitely generated projective
 $R$-modules. Then
\[
\Tor_i^R(P^n, M)=H_i(P^n\otimes B^\cdot)
=H_i(\tilde{B}^\cdot/J^n_{R/k}\tilde{B}^\cdot)
\]
where for $B\in R$-mod, $\tilde{B}$ denotes the 
$R\otimes_k R$-module $R\otimes_k
B$.

Lemma \ref{ArtinRees} below implies that for $i>0$ there exists an $m$ 
such that for all $n$ the maps
\[
H_i(\tilde{B}^\cdot/J^{n+m}_{R/k}\tilde{B}^\cdot)\r 
H_i(\tilde{B}^\cdot/J^{n}_{R/k}\tilde{B}^\cdot)
\]
are zero. This implies that the right hand side of \eqref{CCbasicform}
is zero when $i>0$, and so we are done.
\end{proof}
The following (well-known) lemma, used above, is essentially a more
precise statement of the flatness of completion at an arbitrary ideal
in a Noetherian ring.
\begin{lemmas}
\label{ArtinRees} Let $S$ be a commutative, Noetherian ring  
and let $J\subset S$
be an ideal. Let
\[
P\xrightarrow{\psi} Q\xrightarrow{\phi} T
\]
be an exact sequence of finitely generated $S$-modules. Denote by 
$\psi_n$, $\phi_n$ the induced maps
\[
P/J^nP\xrightarrow{\psi_n} Q/J^nQ\xrightarrow{\phi_n} T/J^n T
\]
and let $M_n$ be the middle 
homology of this sequence. Then there exist $m\ge 0$ 
such that the induced maps $M_{m+n}\r M_n$ are zero for all $n$.
\end{lemmas}

\begin{proof}
This is an immediate consequence
of the Artin-Rees lemma.
The Artin-Rees lemma asserts that given a $J$-adic filtration
on a finitely generated module over a Noetherian ring, 
the induced filtration on any submodule is ``eventually'' $J$-stable
(see any text on commutative algebra, {\it e.g.} \cite{AM}).
In our situation, we see that 
 there exists $m$ such that for all $n$
\[
\phi(Q)\cap J^{n+m}T\subset J^n\phi(Q)
\]
Now let $q\in Q$ represent, modulo $J^{m+n}Q$,  an element of $M_{n+m}$.
Thus 
\[
\phi(q)\in J^{n+m}T\cap \phi(Q)\subset J^n\phi(Q)=\phi(J^n Q)
\]
So there exists $a\in J^n Q$ such that $\phi(q-a)=0$, so that 
 $q-a=\psi(p)$
for some
$p\in P$. This yields
$q\text{ mod }J^nQ=\psi_n(p\text{ mod }J^n P)$ and thus the image of $q$
in $M_n$ is zero.
\end{proof}


Below we indicate some corollaries of Proposition
\ref{CCformgeneral}. They seem to 
suggest that the $R^iD_k(M,N)$ are interesting objects  in their own 
right.
\begin{corollarys}
Let $S \r R$ be a surjective map of  
$k$-algebras satisfying the 
hypotheses of Proposition \ref{CCformgeneral}. Let $M$ and $N$ be $R$-modules,
where $M$ also satisfies the hypotheses of \ref{CCformgeneral}. 
\begin{enumerate} 
\item There is a natural identification
\[
R^iD_{R/k}(M,N)=R^iD_{S/k}({}_S\!M,{}_S\!N)
\]
\item If $R$ is finitely generated over $k$, then
$R^iD_k(M,N)=0$ for $i$ exceeding the minimal number of algebra
generators for $R$ over $k$.
If $R$ is finitely generated over a field,
then $R^iD_k(M,N)=0$ for $i$ exceeding twice the
the Krull dimension of $R$.
\item 
 The  functors $\{R^iD(M, N)\}_i$ form a $\delta$-functor in $M$ 
 on the sub-category of finitely generated $R$ modules which are
 projective as $k$-modules. 
\end{enumerate}
\end{corollarys}

\begin{proof}
The first  assertion  follows because $J_{R/k} = J_{S/k}R$ and in
general for local cohomology $H^i_{JR}(A) = H^i_J({}_SA)$ 
where $J \subset S$ is an ideal in an Noetherian ring,  $R$ is a
Noetherian
$S$-algebra, and $A$ is any $R$ module. 
For (2), note that 
$J$ is generated by
by $x \otimes 1 - 1
\otimes x$ where $x$ ranges through the $k$-algebra generators for
$R$, and that the 
  Krull  dimension of
$R \otimes_k R$ is twice the dimension of $R$, for finitely generated
algebras over a field.
For (3), note that $\Hom_k(-, N)$ is exact, so that
$H^0_J(\Hom_k(-, N))$ is left exact, and thus the 
long exact sequence property follows from that of $H^0_J(-)$.
\end{proof}

\subsection{The Frobenius Functor}
\label{ssec2.1}
We first recall some general facts about 
 rings of prime characteristic $p>0$.
The letter $e$ always denotes a non-negative integer and $q$ always
 denotes an integer power $q = p^e$ of $p$.

For each  $q = p^e$,
the set $R^q$ of $q^{th}$ powers in $R$ forms a subring.
We have a descending chain of subrings of $R$:
$$
R \supset R^p \supset R^{p^2} \supset R^{p^3} \supset \dots\,.
$$ 
Every perfect field contained in $R$ is contained in every subring
$R^q$; in particular, $\ZZ/p\ZZ \subset \cap R^{q}$.
Of course, any $R$ 
module $M$ may be viewed as a module over any one of these subrings $R^q$.

The process of viewing an $R$ module $M$ as a module over
one of these subrings is often described in terms of the 
``Frobenius functor'' on the set of $R$ modules.
The Frobenius map is the ring map $F:R \r R$ sending $r \mapsto r^p$.
We often consider not just the Frobenius map, but also its iterates 
$F^e:R \r R$ sending $r \mapsto r^q$.
Assuming that $R$ is reduced, the Frobenius map $F^e$  is
 an isomorphism onto its image $R^q$.
For any $R$ module $M$, we denote by ${}^e\!M$ the  module $M$ with
its $R$ module structure pulled back via $F^e$. As an abelian group,
${}^e\!M$ is the same as $M$, but its $R$ module structure is given by
$r\cdot m = F^e(r)m = r^qm$.

In this paper, we will frequently be examining the 
structure of $R$ as an $R^q$ module, or equivalently, the structure
of $\,{}^e\!R$ as an $R$ module.
For reduced $R$,
 there is yet a third way to view this 
algebra extension which is sometimes more convenient.
We have an extension $R \subset R^{1/q}$ where $R^{1/q}$ is simply the over-ring
of $q^{th}$ roots of elements in $R$. The Frobenius map $F^e$
affords an isomorphism of $R \hookrightarrow R^{1/q}$ with 
$R^q \hookrightarrow R$.
The descending chain above amounts to an
ascending chain of over-rings of $R$:
$$
 R \subset R^{1/p} \subset R^{1/p^2} \subset \dots.
$$

We use the notation ${}^em$ to indicate that an element $m$ in 
an $R$ module $M$ is being regarded as an element in ${}^e\!M$. In particular,
for an element $x$ in a reduced ring $ R$,
 the element ${}^ex$ in the $R$ module 
${}^eR$ corresponds to the element $x^{1/q}$ in the $R$ module $R^{1/q}$ under
the correspondence discussed above.

Although our formal statements will most often involve 
the Frobenius functor notation, the reader is encouraged to bear all 
three interpretations in mind.
Depending on the context, one  of the following three
equivalent notions may make a particular statement the most transparent:
 the $R^q$ module $R$, the $R$ module ${}^e\!R$, or
the $R$ module $R^{1/q}$.

\subsection{Strong F-Regularity}
\label{strongF}
Throughout this section, we assume
that $R$ is finitely generated as a module over its subring
$R^p$ of $p^{th}$ powers.
This weak assumption (often called ``F-finiteness'')
is satisfied whenever $R$ is a finitely generated algebra over
a perfect field $k$  or whenever $R$ is 
 a complete local Noetherian ring with a perfect residue
field $k$. More generally, $k$ need not even be perfect, as long as 
$[k:k^p] < \infty$. 
F-finiteness  is preserved under localizations.

Two properties of F-finite rings will be of particular
interest to us.

\begin{definitions} \label{Fsplit} If the Frobenius map splits, we say that
$R$ is F-split. That is, $R$ is F-split if   the 
inclusion  $R^p \hookrightarrow R$ splits as a map of $R^p$ modules.
\end{definitions}
This property is also called F-purity in the literature.
Strictly speaking, F-purity, the property that the Frobenius map is pure,
is a weaker condition than F-splitting. However, 
for rings $R$ finitely generated over their subrings $R^p$,
F-splitness is equivalent to  F-purity
(see \cite{HR1}).

Note that if $R^p \subset R$ splits as a map of $R^p$ modules,
then for all $q$, $R^q \subset R$ splits over $R^q$.

A stronger, but closely related property is  {\it
strong F-regularity\/}.
\begin{definitions}
\label{Freg}
A reduced F-finite ring is said to be strongly
F-regular, if for all $c \in R$ not in any minimal prime,
there exists  $q=p^e$ such that the map $R \r R^{1/q}$
sending $1 \mapsto c^{1/q}$ splits as an $R$-module homomorphism.
\end{definitions}

The property of strong F-regularity was introduced by Hochster and Huneke
in \cite{HH2}.
They have conjectured that all ideals of $R$ are tightly 
closed{\footnote{This paper is independent of the theory of tight
    closure, 
but see \cite{HH} for more about this beautiful subject.}}
if and only if $R$ is strongly F-regular. This is known to be true
for Gorenstein rings
\cite{HH2}, for rings of Krull dimension no more than three \cite{W},
and for rings with (at worst) isolated singularities \cite{Mac}.

It is not hard to see that strongly F-regular rings are F-split.
The relationship of this property to the study of differential
operators is apparent from the following fact, proved in \cite{smith}.
\begin{theorems}
\label{dsimple}
Let $R$ be reduced ring finitely generated over $R^p$.
Then 
$R$ is strongly F-regular if and only if 
it is  F-split and it is simple as  a left module over its ring of $\ZZ$ linear
differential operators.
\end{theorems}

Strongly F-regular rings are always Cohen-Macaulay 
and normal \cite{HH2}. In particular, any local 
strongly F-regular ring is a domain.

One more fact  we will need about strong F-regularity follows directly
from the definitions (see  \cite{HH2}):

\begin{theorems}  
\label{SplitF}
If $R \subset S$ splits as a map of $R$ modules and $S$ is strongly F-regular,
then $R$ is also strongly F-regular.
\end{theorems}

\bigskip
\subsection{Rings of Differential Operators in Prime Characteristic}
\label{diffopsp}
Rings of  differential operators
are especially interesting in prime characteristic.
Let $R$ be a reduced ring finitely generated over its subring $R^p$ of
$p^{th}$ powers.
It is not difficult to check that 
\begin{equation}
\label{diffops2}
D_{\mathbb Z}(M, N) = \bigcup_{q = p^e}\Hom_{R^q}(M, N)
\end{equation}
We recall Yekutieli's proof here \cite{Yek} for
convenience of the reader.

Suppose that $\theta$ is a differential operator of order $\le n$.
For each $q>n$, $J_{R/\ZZ}^q \theta = 0$. In particular
 $(r\otimes 1 -  1\otimes r)^q = r^q\otimes 1 -  1\otimes r^q$
kills $\theta$. Thus for all $r \in R$,   $r^q\theta = \theta r^q$, 
and $\theta$ is $R^q$ linear. 

For the converse, suppose that $R$ is generated by
 $\,\xs \,$  as
 an algebra over $R^p$. The same elements generate $R$ over $R^q$,
and thus
for all $e$, $J_{R/R^q}$ is generated by 
$(x_1\otimes 1 - 1 \otimes x_1),\, \dots,\, (x_d\otimes 1 - 1 \otimes x_d)$.
In particular,
$J_{R/R^q}^{dq} \subset 
((x_1^q\otimes 1 - 1 \otimes x_1^q),\, \dots,\, (x_d^q\otimes 1 -
 1 \otimes x_d^q)) = 0$. 
Therefore, any element $\theta \in \Hom_{R^q}(M, N)$ is killed by
$J^{dq}$ and is therefore an $R^q$ linear differential operator.
In particular, $\theta$ is a differential operator from $M$ to $N$ 
linear over any perfect field contained in $R$.

\begin{remarks}
It follows from the proof that 
whenever $R$ is finitely generated over its subring $R^p$ of $p^{th}$ 
powers,
$$
D_{\ZZ}(M, N) = D_{(\ZZ/p\ZZ)}(M, N) = D_{k}(M, N)
= \bigcup_{q = p^e}\Hom_{R^q}(M, N)
$$
where $k$ is any perfect field contained in $R$.  
\end{remarks}

\section{The behavior of modules under Frobenius}
\label{sec3}

\subsection{Finite F-representation type}
\label{ssec31}
Our goal  is to understand the structure of $D_{\ZZ}(R) =
\bigcup_{q = p^e}\End_{R^q}(R).$
We are led to analyze the structure of $R$ as an $R^q$ module,
or equivalently, of ${}^e\!R$ as an $R$ module, as $e \r \infty$.
This task will be  much easier if there are only finitely many 
isomorphism classes of $R$ modules that appear as indecomposable
summands of ${}^e\!R$. 
This observation leads naturally to the introduction of rings
with \emph{Finite F-representation type}.

Finite F-representation type (FFRT) is a representation theoretic
property of commutative rings that makes sense only  in
characteristic $p$. To be  able to state it we need a class of
rings for which the Krull-Schmidt theorem holds. We also need that $R$
is a finitely generated $R^p$-module. Therefore we usually restrict ourselves
to the following two classes of rings.
\begin{itemize}
\item[(A)] Complete local Noetherian rings with  residue field $k$
having the property  $[k:k^p]<\infty$ (the ``complete'' case).
\item[(B)] $\NN$-graded rings of the form $R=k\oplus R_1\oplus
R_2\oplus\cdots$ with $k$ a field such that $[k:k^p]<\infty$ and with
$R$ finitely generated over $k$ (the ``graded'' case).
\end{itemize}

Throughout this section, $R$ will always denote a ring of either type (A)
or type (B).

 Let $M(R)$ stand for the category
of finitely generated $R$-modules 
(resp. finitely generated $\QQ$-graded $R$ modules).
For $M\in M(R)$ we let $[M]$ stand for the isomorphism
class of $M$. 
In the graded case, we do not require the isomorphism to be degree 
preserving:
  $[M]=[N]$ if and only if  $N \cong M(\alpha)$ via a degree preserving isomorphism
 where
$M(\alpha)$ is the graded $R$-module defined by $M(\alpha)_i= M_{i+\alpha}$. If 
$M\in M(R)$ then recall that the notation ${}^e\!M$ means that one is viewing 
$M$ as an $R$-module via restriction of scalars from the Frobenius map.
In the graded case we grade ${}^e\!M$ by putting $({}^e\!M)_\alpha=
M_{p^e\alpha}$.

By the Krull-Schmidt theorem there is a decomposition in $M(R)$
\[
{}^e\!R=M^{(e)}_1\oplus\cdots\oplus M^{(e)}_{n_e}
\]
with $M^{(e)}_i$ indecomposable.
When $R$ is reduced, the reader may find it easier to think about the
(canonically) isomorphic  decomposition of $R^{1/q}$ as an $R$ module where
$q = p^e$. See Section \ref{ssec2.1} 
\begin{definitions}
\label{FFRT}
 We say that $R$ has  Finite F-representation type (FFRT) if the set
\[
\left\{[M^{(e)}_i]\mid e\in \NN, i=1,\ldots,n_e\right\}
\]
is finite.
That is, $R$ has FFRT if there exists a finite set  $\mathcal S$ of 
isomorphism classes of $R$ modules
such that any indecomposable $R$ module summand of $R^{1/q}$, 
for any $q = p^e,\,$ is isomorphic
to some element of $\mathcal S$.
\end{definitions}

If $R$ is regular, then
it is either a polynomial ring or a power series ring over $k$ and 
 ${}^e\!R$ is free over
$R$.
Thus
\begin{observations}
\label{regFFRT}
 Regular rings have finite F-representation type.
\end{observations}

More generally, suppose 
 that $R$ is Cohen-Macaulay. Then ${}^e\! R$ is a Cohen-Macaulay $R$ module
of maximal dimension, 
and hence so are all $M^{(e)}_i$. Now recall that $R$
is said to be of {\em finite representation type} if it has only a
finite number of indecomposable maximal Cohen-Macaulay modules, up to
isomorphism. So we obtain:
\begin{observations}
If $R$ is Cohen-Macaulay and has finite representation type, then $R$
has finite F-representation type.
\end{observations}
This observation applies, for  example, to quadric hypersurfaces in at least
three variables that have an isolated singularity
\cite{Knorrer}.

\bigskip
Below we will show that finite representation type is much stronger
that FFRT. Our first result in this direction is the following.
\begin{propositions} Assume that $R\subset S$ is an inclusion of
rings of type (A) or (B) such that $S$ is a finite
$R$-module and such that ${}^f\!R$ is a direct summand of ${}^f\!S$ as
$R$-modules
for some $f$. Then if $S$ has FFRT,  then so does $R$.
\end{propositions}
\begin{proof}
We have to show that the indecomposable summands of all ${}^e\!R$ are
finite in number as $e \r \infty$. It is sufficient to consider the case $e\ge
f$. For each $e\ge f, $
 ${}^e\!R$ is a direct summand of ${}^e\!S$ as an $R$ module.
 Let $M$ be an
indecomposable $R$ module summand of ${}^e\!R$, $e\ge f$. 
Because ${}^e\!R$ is an $R$ module summand of ${}^e\!S$, 
 $M$ is also an $R$ module
summand of ${}^e\!S$.
Since the $S$ module decomposition of ${}^e\!S$ is automatically an
$R$ module decomposition as well (though not necessarily into indecomposables),
  $M$ must be an $R$ module
 summand of some indecomposable
$S$ summand of ${}^e\!S$. Since the latter are finite in number by
hypothesis, the number of possibilities for $M$ is also finite.
\end{proof}
This proposition shows that rings of invariants  of regular rings under
finite groups of order prime to the characteristic will have FFRT. On
the other hand, in dimension three or higher, these rings of invariants  do
not usually have finite representation type.
Indeed, every ring with finite representation type must have an
isolated singularity \cite{auslander}, whereas rings of invariants need not.
Yoshino's book is  a 
 good source of information about rings of finite representation type
 \cite{Yoshino}.

Now we discuss FFRT in the graded case. The following lemma shows that
there is a  restriction, independent of $e$, on the degrees of
the generators of ${}^e\! R$ as $R$-module.

\begin{lemmas}
\label{degreebound}
Let $R$ be a graded ring as in (B). Assume that $R$ is generated as
$k$-algebra by generators in degrees $a_1,a_2,\ldots,a_n$. Then
${}^e\!R$ is generated as $R$-module by generators whose degrees are
contained in the half open interval $[0,a_1+\cdots+a_n[$.
\end{lemmas}

\begin{proof}
We may assume that $R$ is a polynomial ring $k[X_1,\ldots,X_n]$, $\deg
X_i>0$. Then $R^{1/q}\cong k^{1/q}[X_1^{1/q},\ldots,X_n^{1/q}]$,
$q=p^e$. Let $u_1,\ldots,u_q$ be the generators of $k^{1/q}$ as
$k$-vector space. Then the generators of $R^{1/q}$ as $R$-module are
given by $u_iX_1^{\alpha_1/q}\cdots X_n^{\alpha_n/q}$,
$0\le\alpha_j\le q-1$. The degree of such a generator is given by
$(\sum a_i \alpha_i)/q$, which lies indeed in $[0,a_1+\cdots+a_n[$.
\end{proof}

\begin{propositions}
\label{prop33}
Assume that $R\subset S$ are graded rings satisfying (B). Assume
furthermore that $[S_0:R_0]<\infty$ and that ${}^f\!R$ is an $R$-module
direct summand
 of ${}^f\!S$  for some $f$. Then  if $S$ has FFRT, so
does $R$.
\end{propositions}
\begin{proof}
We have to describe the indecomposable $R$ module summands of 
${}^e\! R$.
 As before,
we may assume $e\ge f$.

Let $N_1,\ldots, N_t$ be the list of indecomposable $S$ module summands of 
${}^e\! S$, up
to isomorphism and shift of degree.
 Without loss of generality, we may, and we will, assume that the grading
on the $N_i$ is normalized in such a way that $(N_i)_0\neq 0$ and
$(N_i)_\alpha= 0$ for $\alpha<0$. 

For any $\QQ$-graded $S$-module $M$ we have a decomposition as $S$-module
\[
M=\bigoplus_{\alpha\in [0,1[\cap \QQ}[M]_{\alpha\text{ mod }\ZZ}
\]
where the index  $\alpha {\text { mod }} \ZZ \in \QQ/\ZZ$ gives rise to  
$[M]_{\alpha\text{ mod }\ZZ}  = \bigoplus_{n\in\ZZ}M_{\alpha+n}$.
Applying this decomposition to the indecomposable $N_i$ we obtain $N_i=
[N_i]_{0\text{ mod }\ZZ}$. Thus the $N_i$ are $\NN$-graded.

Now let $M$ be an indecomposable graded summand of ${}^e\!R$, $e\ge
f$.  Note that any minimal generator of $M$ is a minimal generator
of ${}^e\!R$, so  by  Lemma \ref{degreebound}, the degrees of the 
minimal generators of $M$ are bounded above by $a$, where $a$ is independent
of $e$. Furthermore,  $M$ is also 
an $R$-module summand of ${}^e\!S$ and hence an $R$-module summand of some
$N_i(\alpha)$, $\alpha\le 0$. Thus $M'=M(-\alpha)$ is a $R$-summand of
$N_i$. Note that the  
 minimal generators of  $M'$ have degrees 
bounded by $\le a+\alpha\le a$, where $a$ is independent of $e$. 

So we have shown that any direct $R$-summand of ${}^e\!R$ is, after shifting
degrees, a direct summand of some $N_i$ and is generated in degree $\le a$.
Now let
\begin{equation}
\label{CClist}
M'_1,M'_2,\ldots M'_n,\ldots
\end{equation}
be a (possibly infinite) list of non-isomorphic direct summands of $N_i$,
generated in degree $\le a$. By the Krull-Schmidt theorem, $M_1'\oplus
\cdots\oplus M'_n$
is a direct summand of $N_i$ for all $n$. Now since all ${[N_i]}_k$ are finite
dimensional $S_0$-vector spaces, they also are finite dimensional $R_0$
vector spaces. Furthermore all $M'_i$ contain at least one element in degree
$\le a$. So comparing dimensions in degrees $\le a$, we see that the list
\eqref{CClist} must be finite.

Hence there are only a finite number of non-isomorphic direct summands 
generated in degree $\le a$ of a given $N_i$. Since the $N_i$ are also
finite in number we are done.
\end{proof}

This proposition applies for example to rings of invariants of
{\em linearly} reductive groups. 
A linearly reductive group is one for which all  finite
dimensional representations are completely reducible.
In characteristic $p>0$, the only linearly reductive groups are 
extensions of finite groups whose order is prime to $p$ by tori \cite{Nagata}.
 For groups defined
over a field of characteristic zero, linearly reductive is 
equivalent to reductive \cite{Haboush}.
In the next section, we discuss rings of invariants of linearly reductive
groups in characteristic $p$ in some detail.

\begin{examples} 
The cubic cone, given by $R=k[X,Y,Z]/(X^3+Y^3+Z^3)$, is example of a
graded ring which does not have FFRT, at least when $k$ is
algebraically closed of 
characteristic 
$p = 1 $ mod $3\ZZ$. 
Indeed, in \cite{Tango}, Tango uses Atiyah's classification
of vector bundles to give an explicit decomposition of the structure sheaf
$F_*\mathcal O_X$  of an elliptic curve $X$ viewed as an $\mathcal O_X$ module
via the Frobenius endomorphism.
His method easily generalizes to give a decomposition of $F^e_* \mathcal O_X$
over $\mathcal O_X$, and hence of ${}^e\!R$ over $R$, where
$R$ is the coordinate ring for   $X \subset \PP^2$ ({\em e.g.} the cubic 
cone above).
  In the case where $X$ is not super-singular (equivalently, 
where $R$ is F-split), the decomposition is 
$$
F^e_* \mathcal O_X \cong \bigoplus_{i=0}^{q-1} \mathcal L_i,
$$
where the $\mathcal L_i$ run through the various (non-isomorphic) degree zero 
line bundles corresponding to the $q = p^e$ distinct $q$-torsion 
points\footnote{The $p^e$-torsion points 
on a non-super-singular 
elliptic curve form a cyclic group of  order $p^e$ \cite[p 137]{silverman}.}
 of $X$.
This decomposition corresponds to an $R$ module decomposition
 $({}^e\! R)_{0\text{ mod }\ZZ} \cong \bigoplus_{i=0}^{q-1} M_i$, where the $M_i$ are indecomposable
graded $R$ modules whose sheafifications produce the $\mathcal L_i$.

There are infinitely many non-isomorphic $M_i$ as $e \r \infty$. To
see this, note that if $M_i \cong M_{j}(n)$ for some integer $n$, then
$\mathcal L_i \cong \mathcal L_j \otimes \mathcal O_X(n)$ as $\mathcal
O_X$ modules. Because each $\mathcal L_i$ has degree zero, this is
impossible unless $n = 0$, a contradiction if $i \neq j$. Since
$({}^e\! R)_{0\text{ mod }\ZZ}$ is a direct summand of ${}^e\!R$, we infer that $R$
does not have FFRT.

\end{examples}

\subsection{Invariants under linearly reductive groups}
\label{ssec32}
In this section we
give a more explicit description of the indecomposable summands
of ${}^e\!R$ when $R$ is a ring of invariants for the action 
of a linearly reductive group  on the symmetric algebra $S = S(W)$
of some 
representation $W$ over a perfect field $k$. 
This section is independent of the rest of the paper.

\medskip
Let $U$ be another  representation of $G$. Then $R(U)=(S\otimes U)^G$
is an $R$-module, called a \emph{module of
covariants}.
 Clearly $R(U\oplus V)=R(U)\oplus R(V)$. For the basic properties of modules of
  covariants (in characteristic zero) see 
\cite{VdB9}\cite{VdB3}.
Below we show
that there is a finite set of irreducible representations
$(U_i)_{i=1,\ldots,n}$ such that all ${}^e\!R$ can be written as
direct sums of copies of $R(U_i)$.
These $U_i$  will be described in terms of \emph{Frobenius
twists} of the original representation, so we now digress to review this idea;
see  \cite{jantzen} for a more detailed discussion.

\subsubsection{Frobenius Twisted Representations}
\label{Ftwist}
Let $G$ be a group and let $U$ be any finite dimensional 
representation over a  perfect field $k$.
As before, we may consider the
$k$ module ${}^eU$ defined as the abelian group $U$ but with $k$ structure
defined via Frobenius: $\lambda \cdot v = \lambda^{p^e} v $ for $\lambda \in k$
and $v \in {}^e U$.  Because $G \r GL_k(U) \r GL_k({}^e U)$,
it is clear that ${}^e U$ is also a representation of $G$.
This representation is called a {\em Frobenius Twist} of $U$.
It has traditionally been more standard in representation theory circles
to denote this representation by $U^{(-e)}$. In order to be consistent
 with this tradition, we will adopt this notation here as well.

The apparent anomaly of notation is motivated by the following observation.
Fix a basis $e_1, \dots, e_n$ for $U$. Suppose that the representation
of $G$ is defined by $g \cdot e_i = \sum_{i=1}^n c_{ij}(g) e_j$ for some 
functions $c_{ij}$ on $G$. Because the $e_i$ also form a basis for
${}^eU$
(assuming $k$ is perfect),
 the action of $G$ on ${}^eU$  is described by
$g \cdot e_i = \sum_{i=1}^n (c_{ij}(g))^{p^{-e}} \cdot e_j$.  
Therefore, a collection of matrices defining the representation $U^{(-e)}
(= {}^e U)$ is  given by raising each of the entries in a defining collection
of matrices for $U$ to the $p^{-e}$-th power.

This explicit description of the representation $U^{(-e)}$ brings
to light an important point. Suppose that $G \r GL_k(U)$ is an algebraic
representation, \emph{i.e.,} that 
$G$ is an algebraic group defined over $k$ and that 
the functions $c_{ij}$ above are
rational functions on $G$. 
The twisted representation $G \r GL_k({}^eU)$ is not {\em a priori}
 an algebraic
representation, as the functions $c_{ij}^{p^{-e}}$ need not be in $k[G]$.
The representation $U^{(e)}$ {\em is} algebraic if $e\geq 0$. For
$e < 0$, $U^{(e)}$ is not algebraic in general, though $U^{(e)}$ is easily 
checked algebraic for all $e \in \ZZ$ when
$G$ is a finite group.
 
When $G$ is a torus, the fact that the representations ${}^eU$ are not
algebraic can be side-stepped by appropriately grading the representation $U$.
Thus, though the ideas are similar, there are some slight differences between
our treatments of the finite group case and the torus case. 
 With a little effort, but at
the cost of some extra technicality, similar results can be proved for 
the case where $G$ is an extension of  a    finite group by a torus. We
leave this to the reader.

\subsubsection{$G$ finite}
\label{finitecase}
Let $G$ be a finite group and let $W$ be  a finite dimensional
 $G$ representation over a perfect field $k$. Let $S$ be the 
symmetric algebra $S = S(W)$ and let
 $S_+=\bigoplus_{n>0} S_n$. The notation $S_+^{[p^e]}$ denotes the
ideal of $S$ generated by the $p^e$-th powers of the elements in $S_+$.
\begin{propositions}
\label{prop31}
 With notation as above, set  $R = S^G$, and assume $p$ does not
 divide $|G|$.
\begin{enumerate}
\item
As $R$-modules we have
\[
{}^e\!R \cong \left((S/S_+^{\,[p^e]})^{(-e)}\otimes_k S\right)^G
\]
In particular, ${}^e\!R$ is a module of covariants.
\item
Let $U_1,\ldots, U_n$ be the list of irreducible representations of
$G$ occurring in $(S^lW)^{(f)}$ for some $l,f$. Then ${}^e\!R$ is a
direct sum of modules of covariants  $R(U_i)$ 
and conversely every such module of covariants 
$R(U_i)$ is a direct summand of some ${}^e\!R$.
\end{enumerate}
\end{propositions}
\begin{proof}
For (1), note that we 
 have a $G$-equivariant surjective map $S\rightarrow
S/S^{\,[p^e]}_+$. Since $G$ is linearly reductive this map splits
$G$-equivariantly. This yields a $G$-equivariant map of $S^{p^e}$-modules
\begin{equation}
\label{eq31}
S/S^{\,[p^e]}_+ \otimes_k S^{p^e} \rightarrow S
\end{equation}
Since these are both projective $S^{p^e}$-modules, and \eqref{eq31}
is clearly an isomorphism after the base change
 $\otimes_{S^{p^e}} S^{p^e}/S_+^{p^e}$,
 Nakayama's lemma
yields that \eqref{eq31} is a isomorphism.
Hence
\begin{equation}
\label{eq315}
 S/S^{\,[p^e]}_+ \otimes_k S^{p^e} \cong {}^e\!S
\end{equation}
as $(G,S)$ modules, where $S$ acts on the left hand side of
\eqref{eq315} through
$F^{e}$,
that is, $s\cdot(\bar x \otimes r^{p^e}) = (\bar x \otimes s^{p^e}r^{p^e})$. 
Furthermore,
the bijection 
\begin{equation}
\label{eq32}
(S/S^{\,[p^e]}_+)^{(-e)}\otimes S \r S/S^{\,[p^e]}_+\otimes S^{p^e}:
\bar x \otimes s \mapsto \bar x \otimes s^{p^e}
\end{equation}
is easily checked to be an isomorphism 
of $(G,S)$-modules, where $S$ now acts in the natural way on the
left hand side of \eqref{eq32}.
 Since $({}^e\!S)^G={}^e\!R$, this proves~(1).

For (2), we have to decompose $(S/S^{\,[p^e]})^{(-e)}$ into irreducible
$G$-representations. Now all irreducible representations of $G$ are
defined over some finite extension of the prime field. Hence there
exists a  $u$ such that for every finite dimensional
$G$-representation one has $U^{(u)}\cong U$.

Now every $U_i$ occurs by definition in some $(S^lW)^{(-e)}$. Because
of the remark in the previous paragraph we may assume $e\gg 0$. Then
$U_i$ also occurs in $(S/S_+^{\,[p^e]})^{(-e)}$. Hence $R(U_i)$ is a
direct summand of ${}^e\!R$.

Conversely $(S/S^{\,[p^e]})^{(-e)}$ is a quotient of $S^{(-e)}$ and
hence is a direct sum of the $U_i$. So ${}^e\!R$ is a direct sum of the
$R(U_i)$.\qed
\def\qed{}
\end{proof}

\subsubsection{G a torus} 
Let $G$ be a torus,  split over the perfect ground  field $k$.
 Let $X(G)$ be the character group of $G$; we will use additive notation
in working with $X(G)$.
The notation $X(G)_{\QQ}$ stands for  the group $X(G) \otimes_{\ZZ} \QQ$. 

 For $\chi\in X(G)$ we denote by $L_\chi$ the corresponding
one-dimensional $G$-representation.
Let $W$ be any finite dimensional representation of $G$ over $k$. 
Diagonalizing the action
of $G$ on $W$, we have a decomposition 
 $W \cong \oplus L_{\alpha}$ where $\alpha \in X(G)$
 runs through  the weights of $W$.
The symmetric algebra $S = S(W) = S(\oplus L_{\alpha})$
 is an $X(G)$ graded $k$-algebra.
The graded component $S_{\chi}$ corresponding to $\chi \in X(G)$ consists
of all elements $s \in S$ for which $g \cdot s = \chi(g) s$.
The $k$-algebra ${}^e\!S = S({}^e\!W)$ can be given an $X(G)_{\QQ}$
grading.  In general, for any $X(G)_{\QQ}$ graded object
$U$, let ${}^e\!U$ denote the same abelian group $U$, but with 
the grading shrunk by $[{}^e\!U]_\alpha= [U]_{\alpha p^e}$.

The decomposition of ${}^e\!R$ will be described in terms of the 
\emph{strongly critical characters} of $W$.
Recall the definition.
\begin{definitions}
\label{stcrit}
 Let $\alpha_1,\ldots,\alpha_d$ be
the weights of $W$. A character  $\chi\in X(G)$ is {\em strongly
critical} with respect to $W$ if $\chi=\sum_i u_i \alpha_i$ in
$X(G)_\QQ$ with $u_i\in ]-1,0]\cap \QQ$.
\end{definitions}
 The importance of this
property is that it is a tractable criterion for $R(L_{-\chi})$ to be
Cohen-Macaulay
\cite{St1}\cite{VdB1}.

The following result covers the torus case.
\begin{propositions} Let $G$, $W$ be as above.  
Let $\chi_1,\ldots,\chi_n\in X(G)$ be those characters that are
strongly critical with respect to $W$ and that are weights of some
$S^lW$. Then all ${}^e\!R$ are directs sums of the modules of covariants
 $R(L_{-\chi_i})$ and
conversely every such module of covariants
 $R(L_{-\chi_i})$ is a direct summand of some ${}^e\!R$.
\end{propositions}
\begin{proof}
The proof parallels that of Proposition \ref{prop31}. First diagonalize the 
action
of $G$ on $W$ so as to assume that 
$S=k[x_1,\ldots,x_d]$ where $g\cdot x_i=\alpha_i(g)x_i$ for $g\in
G$.  For any character $\chi \in X(G),$ it is easy to check that 
the module of covariants  $R(L_\chi)$ is $R$-isomorphic to the graded piece 
$S_{-\chi}$.
Note that this is zero unless $-\chi$ is some weight of some $S^l(W)$.
 For an $X(G)_{\QQ}$ graded object $U$,
we employ the
 notation $\Supp U$  to indicate 
the set 
$\{\chi \in X(G)_\QQ\mid U_\chi\neq 0\}$.

The analog of \eqref{eq315}\eqref{eq32} is
\[
{}^e\!S\cong {}^e\!(S/S^{\,[p^e]}_+)\otimes_k S
\]
as $X(G)_\QQ$ graded $S$-modules.
If we decompose $(S/S_+^{[p^e]})$ as $\oplus L_{\chi}$,
then ${}^e(S/S_+^{[p^e]}) \cong \oplus \,{}^e\!L_{\chi} \cong 
\oplus L_{{\chi}/{p^e}}.$

 So
\[
{}^e\!R=({}^e\!S)_0\cong \left({}^e\!(S/S_+^{\,[p^e]})\otimes_k S\right)_0
 \cong
\bigoplus_\chi R(L_{{-\chi}/{p^e}})
\]
where in this direct sum $\frac{\chi}{p^e}$ runs through 
\[
-\Supp\left({}^e\!(S/S_+^{\,[p^e]})\right)\cap \Supp S
\]
with appropriate multiplicities.

Now 
\begin{equation}
\label{eq33}
-\Supp\left({}^e\!(S/S_+^{\,[p^e]})\right)
=
\left\{\sum -\frac{v_i}{p^e}\alpha_i\mid 0\le v_i\le p^e-1\right\}
\end{equation}
and this is contained in the set of strongly critical characters.

Conversely let $\chi=\sum_i u_i\alpha_i$, $u_i\in ]-1,0]$ be a
strongly critical weight, contained in $\Supp S$. Then $\chi=\sum_i
a_i\alpha_i$, $a_i\in\NN$. Assume that $n$ is a common denominator for
the $(u_i)_i$. Then we 
can find for all $\epsilon>0$, $u$ and $e$ such
that 
\[
1\le \frac{un}{p^e}\le 1+\epsilon
\]
Put $\eta=\frac{un}{p^e}-1$ and put $v_i=u_i+\eta(u_i-a_i)$. Clearly
$\chi=\sum_i v_i\alpha_i$ and furthermore $v_i\in ]-1,0]$ if we choose
$\epsilon$ small enough. Since the denominators of the $v_i$ divide
$p^e$ we find that $\chi\in \text{RHS}\eqref{eq33}$. This concludes the
proof.
\end{proof}

\subsection{Growth}
\label{ssec33}
In this section, we say more about the structure of $R$ as 
an $R^q$ module for rings of finite F-representation type. 

Assume that $R$ has finite F-representation type. Then there is a finite set of
indecomposable modules in $M(R)$ such that for each $e$,
\[
{}^e\! R=M_1^{\oplus a_1}\oplus\cdots\oplus M^{\oplus a_n}_n
\]
(with appropriate shifts in the graded case).

Obviously the \emph{multiplicities} $a_1,\ldots, a_n$ depend on
$e$. In this section we show that under the  additional hypothesis of
strong F-regularity these multiplicities
grow like $p^{de}$ where $d$ is the Krull dimension of $R$.

For $M,N\in M(R)$ indecomposable, $e\in\NN$ let us denote by
$m(e,M,N)$ the multiplicity of $M$ in ${}^e\!N$.  Note that
\begin{equation}
\label{eq34}
m(e+f,M,N)=\sum_K m(e,M,K)m(f,K,N)
\end{equation}
where the sum runs through all isomorphism classes of
indecomposable objects in $M(R)$.
This formula follows from the observation that
 $  {}^{e+f}N={}^e({}^fN) = \oplus_K {}^eK^{m(f, K, N)}=
\oplus_K\oplus_M M^{m(e,M,K)m(f,K,N)}$ but also
${}^{e+f}N=\oplus_M M^{m(e+f,M,N)}$.

\begin{propositions}
\label{prop32} Let  $R$ be complete or graded as in (A) or
(B). Assume in addition that $R$ is  strongly F-regular and has
FFRT. Let $d$ be the Krull dimension of $R$ and let
$M_1,\ldots,M_n$ be the list of indecomposable summands of
$\,{}^e\!R$ as $e$ ranges over all natural numbers. Then
\[
\lim_{e\r \infty} \frac{m(e,M_i,M_j)}{p^{de}}
\]
exists and is strictly positive.
\end{propositions}

\begin{proof}
Put $e_{ij}=m(1,M_i,M_j)$ and let $E$ be the 
$n \times n$ matrix whose $(ij)^{th}$ entry is  
$(e_{ij}).$ From \eqref{eq34} it
follows that 
\begin{equation}
\label{eq35}
m(e,M_i,M_j)=(E^e)_{ij},
\end{equation}
the $(ij)^{th}$ entry of 
the matrix $E^e$ obtained be multiplying $E$ by itself $e$ times.
 We claim that some power $E^u$ of $E$ has
strictly positive entries. This means that $m(u,M_i,M_j)>0$ for all $i,j$. 

Since $m(u+v,M_i,M_j)\ge m(u,M_i,R)m(v,R,M_j)$, it is sufficient to
show that there exists $u_0$ such that for all $u \geq u_0$,
and all $i$, we have  $m(u,M_i,R)>0$  and  $m(u,R,M_i)>0$.

Because $R$ is strongly F-regular, it is F-split, and thus 
each 
${}^e\!R$ is a direct summand of ${}^{e+1}\!R$. 
In particular, once an $R$ module $M_i$ appears as a direct summand 
of some ${}^eR$, it appears as a direct summand also of each ${}^fR$ for
$f\ge e$.
Hence, for all  $u\gg 0$,  $m(u,M_i,R)>0$.

Now we consider $m(u,R,M_i)$. As above, $M_i$ is a direct summand of
${}^{e_i}\!R$. Pick an arbitrary $0\neq c\in M_i$. By the definition of 
strongly F-regular the map $R\rightarrow {}^{f_i}\!M_i \hookrightarrow
{}^{e_i+f_i} \!R$ given by sending $1$ to $c$ will split for some
$f_i$. So $R$ is a direct summand of some  ${}^{f_i}\!M_i$.
 Consider any $u \geq
\max \{f_i\}$. Since $R$ is F-split, $R$ is a direct summand of
${}^1\!R$; hence if $R$ is a direct summand of ${}^f\!M$ for some $f$, 
then $R$ is also a  direct summand of ${}^{f+1}\!M$. So $R$ is a
direct summand of all ${}^u\!M_i$ and consequently $m(u,R,M_i)>0$ for all $i$.

Because strongly F-regular rings are normal,
 the fact that $R$  is (graded) local implies that $R$ is a domain. 
So objects in
$M(R)$ have a well defined rank. Since
\[
\rk_R {}^e\!M=\rk_R{}^e\!R\cdot \rk_R M = p^{de} \rk_R M
\]
we obtain the formula
\begin{equation}
\label{eq36}
\sum_i m(1, M_i,M_j)\rk M_i=p^d \rk M_j.
\end{equation}

Let $w$ be the row vector in ${\QQ}^n$ whose $i^{th}$ component is 
the integer $\rk M_i$. With this notation, formula (\ref{eq36}) becomes
$\,wE = p^d w$.\,
Hence $w$ 
is a row eigenvector for $E$
 with strictly positive entries. 
The standard linear algebra lemma below (\ref{lem31}) guarantees that
$p^d$ is the eigenvalue of $E$ with largest absolute
value. Furthermore, Lemma \ref{lem31} also ensures that
$\lim_{e\rightarrow \infty} E^e/p^{de}$
exists and has strictly positive entries.
But by \eqref{eq35}, this means that 
\[
\lim_{e\r \infty} \frac{m(e,M_i,M_j)}{p^{de}}
\]
exists and is positive, so the multiplicities grow like $p^{de}$
and the proof is complete.
\end{proof}

\begin{lemmas}
\label{lem31}
Assume that $E$ is a matrix with non-negative real entries such that some power
has \emph{strictly} positive entries. Then $E$ has a unique eigenvalue
$\lambda$ of largest absolute value. Furthermore $\lambda$ is real and
strictly positive. If $v$ is a (row or column) eigenvector
corresponding to $\lambda$ then $v$ can be chosen to have strictly
positive entries and $v$ is the only eigenvector with this
property (up to scalar multiples).
 Furthermore for every vector $w$ with positive entries one
has 
\begin{equation}
\label{eq37}
\lim_{n\r\infty} \frac{E^nw}{\lambda^n}=av
\end{equation}
for some $a>0$.
\end{lemmas}
\begin{proof}
When  $E$ itself has strictly positive entries this is not difficult to prove 
and in any case is  the well-known
Perron-Frobenius theorem
\cite{gantmacher}.
Assume
$E^u>0$. 
The eigenvalues of $E^u$ are of the form $\lambda^u$
with $\lambda$ an eigenvalue of $E$ and the corresponding eigenvectors
are equal.  
Let $\lambda$ be an eigenvalue of $E$ of largest absolute value and let
$v$ be a corresponding column eigenvector (the case of a row eigenvector is 
similar). Then $\lambda^u$ has largest 
absolute values among the eigenvalues of $E^u$ and hence $\lambda^u$ is
uniquely determined by $E^u$ and is strictly positive. 
Furthermore $v$ is an eigenvector
of $E^u$ corresponding to $\lambda^u$, so, changing signs if necessary, 
$v$ may be assumed to have strictly
positive entries.
Hence $v$ is unique (up to scalar multiple)
among the eigenvectors of $E^u$ (or $E$) 
with this property.
 Then the equation $Ev=\lambda v$ together with the fact that
$E$ has non-negative entries implies that $\lambda$ is real and non-negative.
Using $\lambda=(\lambda^u)^{1/u}$ we deduce that $\lambda$ is unique
and strictly positive.

We
also have for any vector $w$ with non-negative components
\[
\lim_{n\r\infty}\frac{E^{nu+r}w}{\lambda^{nu+r}}=\frac{E^r}{\lambda^r}
\left(\lim_{n\r \infty}\frac{E^{nu}w}{\lambda^{nu}}\right)=\frac{a}{\lambda^r}
E^r v=av
\]
for some $a>0$. Since this $a$ is independent of $r$, \eqref{eq37} follows.
\end{proof}

\section{Rings of differential operators}
\label{sec4}

\subsection{Endomorphism rings} 
\label{ssec41}

In this section, we summarize  a few basic properties of endomorphism
rings needed in the next section. These properties follow easily from the
definitions, and can be found in a standard text, such as \cite{MR}.
 We assume throughout that $R$ is a ring and that $M(R)$ is a category 
of finitely generated $R$ modules for which  the
 Krull-Schmidt theorem
applies. For example, $R$ could be a ring
 of type (A) or of type (B) as defined in Section 3. In the graded case,
the 
adjectives ``graded'' and ``homogeneous''  implicitly modify all 
modules and maps; see \cite{NVO} for 
 basic facts about graded rings.

Suppose that $M \in M(R).$ Let $\Lambda$ denote the ring $\End_R M$.
Recall that if $M$ is indecomposable, then $\Lambda$ is a local ring,
that is, the set of non-units forms an ideal, which is necessarily 
maximal.
(Readers who usually work with commutative rings are reminded that
a local
 non-commutative 
ring has a unique maximal (two-sided) ideal but the converse 
is false in general.)
If $M$ is a direct sum of $n$ copies of the indecomposable $R$ module $N$,
 then $\Lambda$ is isomorphic
to an $n \times n$ matrix ring with 
entries in $\End_R N$, that is $\Lambda \cong  M(n, \End_R N)$. 
This ring also has a unique maximal ideal: the set of 
all matrices in $  M(n, \End_R N)$ with non-units in each entry.

More generally, suppose that 
$$
M=M_1^{\oplus a_1}\oplus \cdots\oplus M_n^{\oplus a_n}
$$
where each $M_i$ is indecomposable and
$M_i\not\cong M_j$ for $i\neq j$. Then $\phi\in \Lambda$ may be written in
block matrix{\footnote{We 
adopt the convention that matrices act on the right, so that
elements of a direct sum are represented by row matrices.}}
 form $(\phi_{ij})_{i,j=1,\ldots n} 
\in \Hom_R(M_i^{a_i}, M_j^{a_j})$; in particular, 
$\phi_{ij}$ is now
itself an $a_i\times a_j$-matrix with entries in
$\Hom_R(M_i,M_j)$. 
The maximal ideals of this ring are in one-to-one correspondence with the
indecomposable modules $M_i$: the $i^{th}$ maximal ideal is the set of
elements with no units in the $i^{th}$ diagonal block $\End_{R} M_{i}^{a_i}$.

In particular,  
an element  $\phi$ is in the intersection of all 
maximal ideals (that is,  $\phi\in\rad \Lambda$) if and only if each of the
 $a_i \times a_i$ matrices 
$(\phi_{ii})_{i=1,\ldots,n}$ consists of  only non-invertible entries.

In other words, 
\begin{equation}
\label{eq41}
\Lambda/\rad\Lambda \cong \prod^n_{i=1} M({a_i}, D_i).
\end{equation}
where $D_i=\End_R M_i/\rad \End_R M_i$.  If we let $k$ denote
the residue field of $R$, then it is clear that  $[D_i:k]<\infty$ and
from the Krull-Schmidt theorem it follows that $D_i$ is a division
algebra. 
The preceding remarks are summarized by the following lemma.

\begin{lemmas}
\label{lem411}
Let $(R, m)$  be a complete local or graded ring with residue field $k$.
Let $M$ be an arbitrary finitely generated $R$ module and set $\Lambda =
\End_R(M)$.
If $M  \cong M_1^{\oplus a_1}\oplus \cdots\oplus M_n^{\oplus a_n}$
is a decomposition of the $R$ module $M$ into indecomposable $R$ modules, 
then  both the simple $\Lambda$ modules and the maximal two-sided ideals
of $\Lambda$ are indexed by the $M_i$.
In particular,
\begin{enumerate}
\item
Using the isomorphism \eqref{eq41} the simple $\Lambda$-modules can be
identified with  $D^{\oplus a_i}_i$.
\item
The maximal two-sided ideals of $\Lambda$ are of the form 
\[
P_i=\{\phi\in\Lambda \mid \phi_{ii}\text{ contains  only 
non-invertible entries}\}
\] 
where $\Lambda$ is identified with the ring of $n$ by $n$   matrices
whose ${(ij)}^{th}$ entry is an $a_{i}$ by $ {a_j}$ matrix $\phi_{ij}$
 with entries in
$\Hom_R(M_i, M_j)$.  
\end{enumerate}
\end{lemmas}

\subsection{Simplicity of Rings of Differential Operators}
\label{ssec42}
The following theorem is the main result of this paper.
\begin{theorems}
\label{th41} Let $R$ be a commutative Noetherian ring of characteristic $p$,
 either complete (type (A)) or graded 
 (type (B)) as defined in Section 3.
If $R$ is strongly F-regular and has finite F-representation type, 
then the ring of differential operators $D(R)$ 
is a simple ring.
\end{theorems}
\begin{proof} 
 We use the formula \eqref{diffops2}. 
\[
D(R)=\varinjlim_{q}\End_{R^q}(R)=\varinjlim_e \End_R({}^e\!R)
\]
Let us define $\Lambda_e=\End_R({}^e\! R)$. We have $\Lambda_0=R$. In
particular, all $\Lambda_e$-modules are automatically $R$-modules.
Given an element $\phi $ in  the direct limit defining $D(R)$, 
we will denote by ${}^e\!\phi $ the corresponding representative in
 $\Lambda_e$. Thus elements ${}^e\!\phi \in \Lambda_e $ 
and ${}^f\!\phi \in \Lambda_f$ determine the same element $\phi \in D(R)$.

\begin{step} Let $0\neq {}^e\!\phi\in\Lambda_e$. Then there exist  $f\ge e$
such that ${}^f\! \phi \not\in\rad \Lambda_{f}$.
\end{step}
\begin{proof}
Assume that ${}^f\!\phi\in\rad \Lambda_{f}$ for all $f\ge e$. 
This means that the left ideal of  $\Lambda_f$ generated by
${}^f\!\phi$ is contained in $\rad \Lambda_f$.
Therefore, by Nakayama's Lemma, we have an inclusion of {\em proper}
 left $\Lambda_f$ submodules of $R$,  
\[
(\Lambda_{f}{}^f\! \phi)R\subset (\rad \Lambda_{f})R \subsetneq R.
\]
In particular, for every $f \ge e$, 
 $(\Lambda_{f}{}^f\! \phi)R$ is a proper $\Lambda_{f}$ (and hence
$R$)-submodule of $R$.

Now, if $\phi$ is not the zero operator, then $\phi \cdot r \neq 0$ 
for some element $r \in R$. Because $R$ is F-regular, it is simple as
a $D(R)$ module ({\em c.f.}
 Theorem  \ref{dsimple}). Thus there exists some $\psi \in D(R)$ 
such that $(\psi \circ  \phi )\cdot r = 1$.
In particular, 
$({}^f\psi \circ  {}^f\phi )\cdot r = 1$ for all $f\gg 0$, so that 
$(\Lambda_f {}^f\! \phi)  R = R$, a contradiction.
\end{proof}

\begin{step} Let ${}^e\!\phi\in \Lambda_e$. 
Then there exists an $f\ge e$ such that
${}^f\! \phi$ is not contained in any maximal two-sided ideal of
$\Lambda_{f}$.
\end{step}
\begin{proof}
Let $M_1,\ldots, M_n$ be the non-isomorphic indecomposable modules
occurring as direct summands of $({}^e\!R)_{e\ge 0}$. By Proposition
\ref{prop32} (or its proof) there exists a $u>0$ such that
$m(u,M_i,M_j)>0$ for all $i,j$.

By Step 1 we may assume that ${}^e\!\phi\not\in\rad\Lambda_e$. So in the
block  matrix representation $({}^e\phi_{ij})_{ij}$, some ${}^e\phi_{ii}$
contains an entry ${}^e\psi\in\End_R(M_i)$ that is  an automorphism
(see the notation and discussion in Section \ref{ssec41}).

Now  we wish to follow what becomes of ${}^e\psi$ when we map further along
in the direct limit system, say to $\Lambda_u$ where $u\ge e$. Then
${}^u\!\psi$ can again be written in block matrix form
$({}^u\!\psi_{kl})_{kl}$. Since ${}^u\!\psi$
 is invertible it is not contained in
any two-sided maximal ideal of $\End_{R}({}^e\!M_i)$. Hence by
\S\ref{ssec41} every block of the form
$\,{}^u\!\psi_{kk}\,$  for $\,k = 1, 2, \dots, n\,$ must 
 contain invertible entries. Now since ${}^u\! \psi_{kk}$
can be identified with
a  sub-matrix of ${}^u\!\phi_{kk}$, we can conclude that each
of the block matrices  ${}^u\!\phi_{kk}$ along the diagonal of ${}^u\!\phi$
 will contain
invertible entries. So again by \S\ref{ssec41}, ${}^u\!\phi$ is not contained
in any maximal two-sided ideal of $\Lambda_{u}$.
\end{proof}

\begin{step}
Now we conclude the proof of Theorem \ref{th41}.
 Assume that $0\neq I\subset D(R)$ is a
non-trivial two-sided ideal. Then there exists an $e$ such that $I\cap
\Lambda_e\neq 0$ (actually since $R$ is a domain, $e=0$). Pick $0\neq
{}^e\!\phi\in I\cap \Lambda_e$. By the previous step there exists an $f\ge e$
such that ${}^f\!\phi$ is not contained in any maximal two-sided ideal of
$\Lambda_{f}$. 
Since ${}^f\!\phi\in I\cap \Lambda_f$ we deduce $I\cap \Lambda_f=\Lambda_f$. 
Therefore $I = D(R)$, and $D(R)$ is simple.
 \qed
\end{step}
\def\qed{}\end{proof}

\begin{corollarys} Let  $R$ be a graded ring
satisfying (B). Assume that $R\subset S$ is a split inclusion of
graded rings where $S$ is a 
(weighted homogeneous) 
polynomial ring, with the extension of residue fields $R_0 \subset S_0$
being finite.
 Then $D(R)$ is simple.
\end{corollarys}

\begin{proof} Polynomial rings always have finite representation type,
so from  Proposition \ref{prop33}, we know that $R$ has finite
F-representation type. On the other hand, $R$ is strongly F-regular by
Theorem \ref{SplitF}.  The corollary is thus an immediate consequence of
Theorem \ref{th41}.
\end{proof}

A result of this generality is not known in characteristic
zero.
Indeed, a characteristic zero analog of this result would
would give an affirmative answer to the question
 of Levasseur and Stafford (Conjecture \ref{con1}). 
 However it is known in the case that $R=S^G$ where $G$ is either finite
 \cite{Kantor} or a torus \cite{MVdB}.

Theorem \ref{th41} also applies to quadric hypersurfaces.
\begin{corollarys}
Suppose that $R = {k[X_1, X_2, \dots, X_N]}/{(Q)}$ where
$k$ is a field of characteristic $p>2$  
and $Q$ is a quadratic form in the $X_i$'s of rank $\ge 3$ (so in particular
$N\ge 3$).
Then the ring of differential
operators $D_{k}(R)$ is simple.
\end{corollarys}
\begin{proof}
Since rings of differential operators are compatible with flat base change
(see \S\ref{ssec51}) we may assume that $k$ is algebraically closed. We 
may then change coordinates so as to 
 assume that $Q=\sum_{i=1}^m X_i^2$ with $\rk Q\ge 3$.
It is an immediate application of Theorem 4.1 of \cite{S2}
that
 $R$ is strongly
F-regular in all characteristics $p>2$.
(For $p\gg 0$, this was first proven by Fedder \cite{F}.)

It is easy to see
 that the FFRT property is
preserved under taking polynomial extensions, whence  we may assume that $Q$
has maximal rank.
In this case, $R$ has finite 
representation type \cite{Knorrer}, so it certainly has FFRT. 
We can therefore  use Theorem \ref{th41} to conclude that $D(R)$ is simple.
\end{proof} 

 The fact that $D(R)$ is simple  in the  quadric  hypersurface case 
when $k$ is 
characteristic zero was proved in
\cite{Lev2} (see also \cite{Lev1}). The characteristic zero proof uses
the structure of primitive ideals in enveloping algebras 
and is much more intricate
than our proof of the characteristic  $p$-case. 
\subsection{Finite dimensional representations}
\label{ssec43}
Let $R$ be either graded or complete as in (A) or (B). In both case
$R$ contains a copy of its residue field $k$. A $D(R)$-module is said
to be a {\em finite dimensional representation} if it has finite length
as an $R$ module, {\em ie}, if it is finite dimensional as a vector space
over $k$.
 As is to be expected, the finite dimensional
representations of $D(R)$ are   controlled by the direct summands
of ${}^e\!R$. In fact,  in this section we give a necessary and
sufficient, representation theoretic, criterion for $D(R)$ to have no
finite dimensional representations.

Let $W$ be a $D(R)$ module.
The annihilator in $D(R)$ is clearly a two-sided ideal of $D(R)$.
If $W$ 
 is a finite dimensional $D(R)$ module, then because $W$ is in fact killed
by some power of the maximal ideal of $R$, this annihilator 
must be non-trivial. Therefore,
 if $D(R)$ is simple, then it has no finite dimensional
representations (at least if $\dim R>0$). Unfortunately there seems
to be, a priori, no reason why the converse should hold.

For an indecomposable finitely generated  $R$ module $M$, 
 let us define $d(M)$ to be the dimension of the division algebra
$\End(M)/\rad \End(M)$. Of course if $k$ is algebraically closed then
$d(M)$=1.
\begin{propositions}
Let $R$ be complete or graded as in (A) or (B). Then the minimal
dimension of a finite dimensional $D(R)$-representation is given by
\[
u=\sup_e \,\,\{\,\min_{M|{}^e\!R} \,\,\, d(M)m(e,M,R)\,\}
\]
where the $M$ range through all the indecomposable $R$ modules appearing
as direct summands in ${}^e\!R$ with non-zero multiplicity.
If particular $D(R)$ has no finite dimensional representations if and
only if $u=\infty$. The ``sup'' here means supremum.
\end{propositions}
\begin{proof}
This follows from the formula
$
D(R)=\varinjlim \End_R({}^e\!R)
$.
By Lemma \ref{lem411},
the simple $\Lambda_e$ modules are all of the form $D^{\oplus m(e, M, R)}$
where $M$ is one of the indecomposable summands appearing in an $R$ module
decomposition of ${}^e\!R$,  $D$ is the division algebra $\End M/\rad \End M$,
and $m(e, M, R)$ is the multiplicity of $M$ in ${}^e\!R$.
Given a finite dimensional $D(R)$ module, we can, of course, view it as a
finite dimensional $\Lambda_e$ module. Note that
\[
\min_{M; m(e,M,R)\neq 0} d(M)m(e,M,R)
\]
is the minimal dimension of a finite dimensional
$\End_R({}^e\!R)$-module. We can then invoke
Lemma \ref{lem432}. 
\end{proof}

\begin{lemmas}
\label{lem432}
 Assume that $(\Lambda_e)_{e\in\NN}$ is a directed
system of rings containing a field $k$
 each  having only a finite number of non-isomorphic
simple representations. Let $u_e$ be the minimal dimension of a finite
dimensional (over $k$)  
$\Lambda_e$-representation and let  $\Lambda=\varinjlim
\Lambda_e$. Then $u=\sup u_e$ is the minimal dimension of a finite
dimensional $\Lambda$-representation.
\end{lemmas}
\begin{proof}
We first observe that $(u_e)_e$ is a  non-decreasing sequence of integers,
since restriction of 
scalars defines
a dimension preserving functor from $\Lambda_{e}$-mod to $\Lambda_{e-1}$-mod.

\noindent {\bf The case {\boldmath $u=\infty$} } If $W$ is a finite
dimensional $\Lambda$-module then $W$ is also a $\Lambda_e$-module for every 
$e$. So
$\dim W\ge u_e$.
Since this holds for all $e$, $\dim W\ge u=\infty$. Hence in this case
one certainly has $\dim W=u$

\noindent {\bf The case {\boldmath $u<\infty$} }
In this case there exists an $f$ such that for $e\ge f$ one has $u_e=u$.
Define for all $e\ge f$:
\[
\Fscr_e=\{\text{isomorphism classes of simple $\Lambda_e$-modules of
dimension $u$}\}
\]
By hypothesis $\Fscr_e$ is a non-empty finite set. Furthermore the
restriction functor from $\Lambda_e$-mod to $\Lambda_{e-1}$-mod defines a
map $\Fscr_e \rightarrow \Fscr_{e-1}$. Define
$\Fscr=\varprojlim\Fscr_e$. 
But an inverse limit of finite sets is always non-empty, so  
$\Fscr\neq
\emptyset$. 

Let $(W_e)_e\in\Fscr$. 
We may assume that $W_e$ is a simple $\Lambda_e$ module of dimension $u$
representing an isomorphism class in $\Fscr_e$.  Restriction of scalars
(possibly composed with another isomorphism) gives an isomorphism of 
$\Lambda_{e-1}$ modules $W_e \rightarrow W_{e-1}$. Therefore, we have
 isomorphisms $W_e \rightarrow W_{e+1}$  of $\Lambda_e$ modules
for all $e \in \NN$. The direct limit 
 $W=\varinjlim W_i$ is a finite dimensional
$\Lambda$-module of dimension $u$. 

It is clear that there can be no $\Lambda$ module of smaller dimension.
For if $W$ is a $\Lambda$ module of dimension $<u$, then by restriction
of scalars, $W$ is  a $\Lambda_e$ module  of dimension $<u$ 
for each $e$, contrary 
to the definition of $u$.
\end{proof}

\subsection{D-simplicity}
\label{sec5}

Closely related to the simplicity of the ring $D_k(R)$ is the simplicity
of $R$ as a module over $D_k(R)$. If $D_k(R)$ is a simple ring,
then $R$ is simple as a $D_k(R)$ module, as the annihilator of
 any $D_k(R)$ module of the form
$R/I$ is a non-zero two-sided ideal of $D_k(R)$. 
The converse, however, is false; see the example
of Chamarie, Levasseur and Stafford in the introduction of \cite{LSt}.

 How does one verify that a given $k$-algebra $R$ is simple as a 
$D_k(R)$ module?  Clearly $R$ is D-simple if and only if
each  non-zero element $c \in R$ generates all of $R$ as a $D_k(R)$ module;
that is, if and only if for each non-zero $c \in R,$
         there exists some differential
operator $\theta \in D_k(R)$ sending $c$ to  $1 \in R$.
In practice we would like to be able to check this condition,
not for all $c \in R$, but for a \emph{single element} $c \in R$.
We prove in this section that this is indeed the case for a large
class of rings of characteristic $p$, describing a specific $c$ that
governs the $D$-simplicity of $R$.

Note that if $k$ is not a field, but merely a commutative ring, then
$R$ is virtually
 never a simple $D_k(R)$ module, because any ideal of $k$ expands to 
an ideal of $R$ that is stable under $D_k(R)$.  This can be remedied 
by introducing a concept of {\em relative D-simplicity}.
A $k$-algebra $R$ is said to be relatively $D_k$ simple if   
  every non-zero $D_k(R)$ submodule of $R$ has
non-zero intersection with the set of non-zero-divisors in 
$\im k \subset R$.
If $k$ is a field, this is equivalent to  $R$ being a simple $D_k(R)$ 
module.

The next proposition lets us reduce the problem of checking 
relative $D_k$ simplicity to checking that certain elements
can be sent to $k$. We abuse notation throughout by speaking of
elements of $k$ as if $k$ were a subring of $R$, though of course
the structure map $k \rightarrow R$ need not be injective.

\begin{propositions}
\label{pdtest}
Let $R$ be a 
 $k$-algebra, where $k$ is an arbitrary (commutative) ground ring.
Let $c$ be any non-zero-divisor of $R$ such that $R_c$ is
a relatively simple $D_k(R_c)$ module.
Then
 $R$ is relatively 
 simple as a $D_k(R)$ module if and only if for each integer
$n$, there is some element $\theta_n \in D_k(R)$
 such that $\theta_n \cdot c^n$ is a non-zero-divisor in $k$.
\end{propositions}

\begin{proof}
To prove that $R$ is relatively simple as a $D_k(R)$ module, we need only show
that for each non-zero $x \in R$, there is an operator 
$\theta$ sending $x$ to a non-zero-divisor in  $k$.
Let $\frac{x}{1}$ denote the image of $x$ in $R_c$.
Because $c$ is a non-zero-divisor,
$\frac{x}{1}$ is non-zero provided that $x$ is non-zero. 
Because $R_c$ is relatively simple as a $D_k(R_c)$ module, 
there is a differential operator $\theta'$ in $D(R_c)$ taking 
$\frac{x}{1}$ to $\lambda \in k$, where $\lambda $ is not a zero-divisor.
 Since $\,D_k(R_c) \cong D_k(R) \otimes_R R_c,\,$ 
we may assume that $\theta'  = \frac{ \theta}{c^n}$ for some $\theta 
\in D_k(R)$.  Therefore $\theta \cdot x = \lambda c^n$ in $R$.

We have assumed that there is an operator $\theta_n \in D_k(R)$ 
such that $\theta_n \cdot c^n$ is sent to a non-zero-divisor in 
 $k$. Therefore, the composition
$\theta_n \circ \theta \in D_k(R)$ is a differential operator
sending $x$ to a non-zero-divisor in  $k$. The proof  is complete.
\end{proof}

We would like to be able to check that $R$ is simple (or relatively simple)
 over $D_{k}(R)$ by
checking just one condition, not infinitely many. 
We pose the question:
\begin{questions}
Does  there exist an element $c \in R$ such that 
$R$ is simple as a $D_k(R)$ module if and only if $\theta \cdot c = 1$
for some $\theta \in D(R)$?
\end{questions} 
Amazingly, this turns out to be true in characteristic $p>0$.

\begin{theorems}
\label{Dtest}
Assume that $R$ is a domain of characteristic $p>0$ 
 finitely generated over its subring $R^p$.  
Suppose that 
$R$ is module
finite over some F-finite
 regular domain, $T$, such that
the corresponding extension of fraction fields is separable.
Let $c$ be a discriminant of  $R$ over $T$, {\it i.e.\/} 
$c = \det(\trace (r_i r_j))$, where 
$r_1, r_2, \dots, r_n\in R$ is a  basis for $R \otimes_T K$ over $K$, the fraction
field of $T$. Then $R$ is simple as a $D(R)$ module
if and only if there is a differential operator sending $c^2$ to 1.

More generally, $R$ need not be a domain if it is reduced and if it is module 
finite, torsion-free, and  generically smooth over the regular subring $T$.
\end{theorems}

Note that the restrictions on $R$ above are quite weak in general.
By a variant of Noether normalization, any domain finitely generated
over a perfect field satisfies the hypothesis. Similarly,
any complete local domain with a  perfect residue field
 satisfies the hypothesis.

\begin{proof}
Localizing at the discriminant $c$, the map $T_c \hookrightarrow R_c$
becomes \'etale, so $R_c$ is regular.
Because regular rings are $D_{\ZZ/p\ZZ}$ simple, 
 $R_c$ is simple as a $D(R_c)$ module. By Proposition \ref{pdtest},
it therefore suffices to check that there are differential 
operators taking each power $c^n$ to 1. Because multiplication by
$c^{q-n}$ is a differential operator on $R$, there is no
harm is assuming that $n$ is a large power of $p$.

An important property of the discriminant, observed by Hochster and Huneke,
 is that for each power $q$ of $p$,
$cR^{1/q} \subset T^{1/q} \otimes_T R$ (Lemma 6.5 \cite{HH}; see also 
first sentence of proof of Lemma 6.4 in \cite{HH}).
Note also that because $T^{1/q}$ is free over $T$, the ring 
$T^{1/q} \otimes_T R $ is free over $R$. Thus there is an $R$ linear map
$T^{1/q} \otimes_T R \rightarrow R$ which sends $1 \otimes 1$ to 1.
Pre-composing with multiplication by $c$, we 
have an $R$ linear  map
$R^{1/q} \xrightarrow{\text{mult by $c$}} T^{1/q} \otimes_T R 
\rightarrow R$ which sends $1$ to $c$.
Raising everything to the $q^{\text{th}}$ power,
 there is an $R^q$ linear map $\pi$ from $R$
to $R^q$ sending $1$ to $c^q$.

We wish to find a differential operator in $D(R)$ sending $c^q$ to 1.
Suppose that $\theta \in \End_{R^q}{R} \subset D(R)$ sends $c^2$ to 1.
Let $\theta^{[q]} \in \End_{R^{q^2}}{R^q}$ be the operator
$$
R^q \rightarrow R^q
\,\,{\text { sending }}\,\,
x^q \mapsto (\theta \cdot x)^q.
$$
Consider the composition
$$
R \overset{\pi}\to R^q  \overset{\theta^{[q]}}\hookrightarrow R.
$$
The element $c^q$ is sent to $\theta^{[q]} \cdot c^{2q}
 = (\theta \cdot c^2)^q = 1$.
Because this map is $R^{q^2}$ linear, it is a differential operator, and
the proof is complete.
\end{proof}

Let $R$ be a finitely generated algebra over a perfect field $k$
of characteristic $p>0$.
The above theorem tells us that
 there is a \emph{single}
 element, namely
the square of any discriminant of a (separable)
 Noether normalization, which governs 
whether or not $R$ is simple as an $D_k(R)$ module.
We do not know whether or not a similar result holds in characteristic zero
or in the relative setting. An interesting question:
  Why does the discriminant play this special
role with respect to 
differential operators?

\section{Characteristic zero}
\label{sec6}

Ultimately, we would like to be able to prove Conjecture \ref{con1}
in characteristic zero as well.
We would like to accomplish this by 
a  reduction
mod $p$ argument and then invoking
 Theorem \ref{th41}.  In order to do this, three major 
issues need to be addressed:
\begin{enumerate}
\item
An understanding of how differential operators behave mod $p$, at least
for invariant rings;
\item An understanding of  the FFRT property for finitely generated
algebras over a field of characteristic zero;
\item An understanding of strong F-regularity in characteristic zero.
\end{enumerate}

The idea is as follows. Suppose we are given a finitely generated 
$k$-algebra $R \cong {k[X_1, \dots, X_n]}/{(F_1, \dots, F_r)}$
 where $k$ is a field of characteristic zero. We build a finitely generated
$\ZZ$-algebra $A$ that contains all of the elements of $k$ necessary
to define $R$; \emph{i.e.} $A$ should contain all the   coefficients of the 
polynomials $F_1, \dots F_r$ defining the ideal of relations on the generators
for $R$.  We thus have an $A$-algebra 
$R_A ={A[X_1, \dots, X_n]}/{(F_1, \dots, F_r)}$ 
such that the natural
map $R_A \otimes_A k \rightarrow R$ is an isomorphism. By the lemma 
of generic freeness, we may invert a single element of $A$ so as to
 assume that $R_A$ is $A$ free. 
Each of the closed fibers of the map $A \rightarrow R_A$ is a 
finitely generated algebra over a perfect (finite!) field. We think of 
the family of closed fibers as being a model for the original
$k$-algebra $R$.

Several natural questions come to mind:
How does $D_A(R_A) \otimes L$ compare to $D_L(R_A \otimes_A L)$ where
$L$ is an $A$-algebra? If $R$ is simple as a $D_k(R)$-algebra, does
this mean that the closed fibers $R_A \otimes_A A/\mu$ are simple
as $D_{A/\mu}(R_A \otimes_A A/\mu)$ modules? If $R$ is a graded 
 direct summand of
a polynomial ring, is it true that $R_A \otimes_A L$ has finite 
F-representation type where $L$ is a perfect field of characteristic $p$
to which $A$ maps?

\subsection{Reduction to characteristic  $p$ and 
rings of differential operators}
\label{ssec51}

We discuss item (1) in somewhat
more detail. 
Unfortunately, differential operators are not well behaved ``mod $p$''
in general. In \cite{smith}, an example is given to show that in some sense
there are ``more'' differential operators in characteristic $p>0$.
It is the cubic cone ${k[X, Y, Z]}/{(X^3 + Y^3 + Z^3)}$.
Nevertheless, reduction \emph{does} work for some nice classes of rings
such as regular rings and rings of invariants for finite groups. So it
is not unreasonable to expect that reduction should work for some
other good rings, such as rings of invariants under
reductive groups.
 Unfortunately, we have not been able to prove this.

        What we can do, however, is describe the major obstruction against
reduction mod $p$. It is based upon the first derived functor of the
left exact functor $D(R,-)$ (see 
\S\ref{ssec215}). We
write $R^1D(R)$ for $R^1D(R,R)$.

Let $R_A$ be flat and 
finitely generated over $A$. For simplicity we assume that $A$
is a Dedekind domain. Recall that 
\begin{equation}
\label{myformula}
D_A(R_A) = \varinjlim_n\Hom_{R_A}(P^n_{R_A/A}, R_A) 
\end{equation}
where  ${P^n_{R_A/A} = (R_A \otimes_A R_A)/J_{A}^{n+1}}$
with $J_A$ is the kernel of the  
multiplication map $R_A \otimes_A R_A \rightarrow R_A$;
see \ref{diffops}.
Let $L$ be any $A$-algebra.
Because $R_A$ is $A$-flat, it is clear that $J_A \otimes_A L$ 
is the kernel of the corresponding multiplication map for $R_A \otimes_A L$.
Furthermore, by inverting a single element
of $A$, 
 the modules $P^n_{R_A/A} $ for all $n$ may be  assumed 
to be all $A$-flat.
To see this, 
consider the 
short exact sequence 
$0 \rightarrow J^n/J^{n+1} \rightarrow P^{n+1} \rightarrow P^n \rightarrow 0$.
The $A$ modules $J^n/J^{n+1}$ can be assumed $A$ free because they are the
graded pieces of the finitely generated $A$-algebra 
$\Gr_{J} S = S/J \oplus J/J^2 \oplus J^2/J^3  \oplus \dots, $
with $S = R_A \otimes_A R_A$, and this graded ring may be assumed $A$-flat
after inverting a single element of $A$, by the lemma of generic 
freeness. That all the $P^n$ can be assumed $A$-flat now follows immediately
by induction on $n$.
We conclude that 
$P^n_{R_A/A} \otimes_A L = P^n_{R_A \otimes_A L/L}$.

Because tensor product commutes with direct limits, we see that
we therefore have a natural map for any $A$-algebra $L$
\[
D_A(R_A) \otimes_A L \rightarrow D_L(R_A \otimes_A L).
\]
The question is to determine when this map is an isomorphism.
It always is when $L$ is a flat $A$-algebra, for instance, when 
$L$ is the fraction field of $A$. But our main concern is when $L = A/\mu$,
for some $\mu \in \max\Spec A$. On a Zariski dense open subset
of $\max\Spec A$ the $R$-modules $P^N_{R_A/A}$ are $A$-flat.
One can then easily show  using \eqref{myformula}  
and the Universal Coefficient Theorem  
that for the same open set of $\mu$'s in $A$, there is a  short exact sequence
\[
0\r D_A(R_A)\otimes_A A/{\mu} \r D_{A/\mu}
(R_A\otimes_A A/\mu) \r \Tor_1^A(A/\mu,R^1D_A(R_A)) \r
0
\]
So for reduction to work we should have
$\Tor_1^A(A/\mu,R^1D_A(R_A))=0$ on a dense set of maximal ideals $\mu$ in $A$.
 Unfortunately $R^1D(R)$ seems to be very hard
to compute explicitly. One could naively hope that  $R^1D(R)$ is
always zero but
this is contradicted by the quadric hypersurface, which we discuss
below.

\begin{examples}
\label{example51}
Assume that $R=S/I$ where $S$ is a graded polynomial ring over
$A$. 
Using the results in Section \ref{ssec215} we find that
\begin{align*}
R^iD(R)=R^iD({}_S\!R,{}_SR)&=\Ext^i_{D(S)}(R\otimes_S D(S),R\otimes D(S))\\
&=
\Ext^i_{D(S)}(D(S)/ID(S),D(S)/ID(S))
\end{align*}
Let us now assume that $R=A[x_1,\ldots, x_n]/(f)$ where
$f=\sum_ix_i^2$. Then
\[
R^1D(R)=\Ext^1_{D(S)}(D(S)/fD(S) , D(S)/fD(S) )=\frac{D(S)}{fD(S)+D(S)f}
\]
Assume now that $A$ is of characteristic zero and let $k$ be the
algebraic closure of the fraction field of $A$.
 Put $\bar{S}=k\otimes_A S$. We will show
that
\[
k\otimes_A R^1D(R)=\frac{D(\bar{S})}{fD(\bar{S})+D(\bar{S})f}
\]
is not zero.

Since we are in characteristic  zero
\[
D(\bar{S})=k[x_1,\ldots,x_n,\partial_1,\ldots,\partial_n]
\]
where $\partial_i=\frac{\partial\ }{\partial x_i}$. 

Let $O(n)$ be the orthogonal group over $k$, acting in the standard
way on $\bar{S}$. It is easy to see that 
\[
D(\bar{S})^{O(n)}=k[f,g,h]
\]
where
$f=\sum_i x_i^2$,
$g=\sum_i x_i\partial_i$ and
$h=\sum_i \partial_i^2$.
$f,g,h$ satisfy the relations
\begin{align*}
[g,f]&=2f\\
[g,h]&=-2h\\
[h,f]&=4g+2n
\end{align*}
So we see that $\mathfrak{g}=kf+k(g+\frac{n}{2})+kh$ is a Lie algebra
(isomorphic to $\mathfrak{s}\mathfrak{l}_2$). Hence we have a surjective map
$U(\mathfrak{g})\r D(\overline{S})^{O(n)}$. By comparing the dimensions of
the associated graded rings on both sides one sees that there can be
no extra relations among $f,g,h$. So  $D(\overline{S})^{O(n)}$ is
actually isomorphic to the enveloping algebra of $\mathfrak{g}$. Using
this isomorphism we obtain
\begin{equation}
\label{eq51}
(k\otimes_A
R^1D(R))^{O(n)}=\frac{U(\mathfrak{g})}{fU(\mathfrak{g})+U(\mathfrak{g})f}
\end{equation}
Since $f$ is in the augmentation ideal of $U(\mathfrak{g})$, the right
hand side of \eqref{eq51} cannot be zero.
\end{examples}

We pose the following problem.
\begin{questions}
Let $A$ be a domain finitely generated as an algebra over $\ZZ$.
Suppose that $R_A$ is a finitely generated $A$-algebra such that 
$R_A \otimes_A K$ is a direct summand of a regular ring for $K$ the 
fraction field of $A$. Is the  $A$-module 
$\Tor^A_1(A/\mu,R^1D_A(R_A))$ zero for all maximal  $\mu \in \Spec A$
except on some
 Zariski closed subset? What if $R_A \otimes_A K$
is a ring of invariants for a linear  action of a reductive group on
a polynomial ring?
\end{questions}
An affirmative answer would tell us that the ring of differential 
operators on a ring of invariants can be studied by reduction to 
characteristic $p>0$.

\subsection{Reduction mod $p$ and strong F-regularity}
\label{strongred}
Strong F-regularity, though defined in \S\ref{strongF}
as a characteristic $p$ notion, can be made meaningful for 
finitely generated algebras over a field of characteristic zero.

Let $R = {k[X_1, X_2, \dots, X_n]}/{(F_1, F_2, \dots, F_r)}$
be a finitely generated algebra over a field $k$ of characteristic zero.
Choose a finitely generated $\ZZ$ subalgebra of $k$ over which $R$ is 
defined and set
 $R_A =
{A[X_1, X_2, \dots, X_n]}/{(F_1, F_2, \dots, F_r)}$.

We say that $R$ is strongly
F-regular{\footnote{In the tight closure literature, this would  be called  
strongly F-regular type.}}
 if there exists such a choice of $A$
for which all the closed fibers $R_{\bar A}$ 
are strongly F-regular on some non-empty Zariski open set of $\Spec A$.
Similarly, we say that $R$ has finite F-representation type or that 
$R$ is F-split if either of these properties holds on a non-empty 
Zariski open set of $\Spec A$.

We  investigate what sort of $k$-algebras will be
strongly F-regular.
We first note that D-simplicity in characteristic zero implies it on a 
Zariski open set of closed fibers.

\begin{theorems}
\label{Dsimplered}
Let $R_A$ be a reduced
finitely generated $A$-algebra where $A$ is a  domain finitely
generated as a $\ZZ$-algebra. If the $K$-algebra
$R_A \otimes_A K$  (where $K$ is any field containing $ A $) is simple
as a $D_K(R)$ module, then the same is true for almost all closed fibers:
for all maximal $\mu \in \Spec{A}$ in a Zariski dense open set,
the ring $\bar R$ is simple as a $D_{\bar A}(\bar R)$ module, where the
``bar'' indicates reduction mod $\mu$.
\end{theorems}

\begin{proof}
We first note that it is sufficient to prove the theorem when 
$K$ is the fraction field of $A$. Let $L$ be the fraction field of $A$.
Because  $L \subset K$ splits over $L$, the inclusion 
 $D_L(R_A \otimes L)  \cong D_A(R_A) \otimes L \subset 
D_K(R_A \otimes K) \cong D_A(R_A) \otimes K$
splits over $D_L(R_A \otimes L)$.
Assume that $R_A\otimes K$ is $D_K(R)$ simple.
Then for each element $c \otimes 1 \in  R_A \otimes L  \subset R_A \otimes K$,
we can find an operator in $D_A(R_A) \otimes K$ sending $c \otimes 1 $ to 1.
This operator is of the form $\sum \theta_i \otimes \kappa_i$ where
$\theta_i \in D_A(R_A)$ and $\kappa_i \in K$. 
Applying the splitting above, there is a corresponding operator
in $D_L(R_A \otimes L)$, $\sum \theta_i \otimes \lambda_i$ where 
$\lambda_i \in L$. Using the splitting we see that
 $\sum \kappa_i \theta_i \cdot c = \sum
\lambda_i \theta_i \cdot c = 1$,
so that $R_A \otimes L$ is simple over $D_L(R_A \otimes L)$.
Thus we may assume that $K$ is the fraction field of $A$.

The point is that 
(after possibly inverting a single element of $A$) 
we can find a Noether normalization for $R_A$, 
{\it i.e.\/}, a polynomial  subring $T_A \hookrightarrow R_A$
over which  $R_A$ is module finite and generically smooth.
The discriminant of this algebra extension is an element $c_A$
of $T_A$ whose image in  almost all 
 the fibers is the corresponding discriminant for
the algebra extension in the fibers (see \cite{HH3}, Discussion 2.4.5).  
Thus, 
by Theorem \ref{Dtest},
 in order to check that $\bar R$ is simple over $D(\bar R)$,
we need only find a differential operator in $D(R)$ taking $\bar c_A^2$
to 1. 

Because the generic fiber is D-simple, there is a differential operator
$\theta \in D(R \otimes K)$  sending $c_A$ to 1. 
Without loss of generality, we may assume that 
$\theta \in D_A(R_A)$. But because for almost all fibers we have 
$D_A(R_A) \otimes \bar A \subset  D_{\bar A}(\bar R)$, this means that 
$\bar \theta $ is differential operator on $\bar R$ sending $\bar c^2$ to 1,
as needed.
\end{proof}

Strong F-regularity in characteristic zero is decidedly more subtle.
It is no longer clear that  a finitely generated $k$-algebra that is 
a direct summand of a regular ring in 
characteristic zero will be strongly F-regular.
Although the preceding result shows that the D-simplicity descends to
characteristic $p$, the issue of F-splitting is harder to deal with.
Using  results of Hochster and Roberts \cite{HR1},
one can deduce the F-splitting for graded Gorenstein direct summands of
a polynomial ring, at least on a dense set of closed fibers.
Using the more recent theory of tight closure in characteristic zero
due to Hochster and Huneke \cite{HH3}, one can prove strong F-regularity
for a Zariski open set of fibers. We record this proof below.
\begin{lemmas}
\label{gorcase}
Let $R$ be a 
graded Gorenstein ring 
finitely generated as a  $k$-algebra, where $k$ is
a field of characteristic zero. Assume that $R$ is a direct summand
(as a graded $R$ module) of a graded regular ring.
Then $R$ is strongly F-regular on a Zariski open set of fibers of any finitely
generated $\ZZ$-algebra  $A$
over which $R$ is defined.
\end{lemmas}

\begin{proof}
By Proposition 3.1 of \cite{smith},
 $R$ will be simple as a $D_k(R)$ module. By Proposition \ref{Dsimplered},
it follows that for any choice of $A$, the closed fibers 
$R_A \otimes_A {A/\mu}$ are  simple as $D_{A/\mu}(R_A \otimes_A A/\mu)$
modules on a Zariski dense open set of $\max\Spec A$. 
 Therefore, by Theorem \ref{dsimple}, it is enough to show that
a non-empty Zariski open set of the fibers is F-split.

We choose to do this using the theory of tight closure in characteristic zero
(see \cite{HH3}). The point is that because $R$ is pure in a regular ring,
all ideals of $R$ are tightly closed (in the characteristic zero theory). 
Unfortunately, this is not known to imply in general 
that  a Zariski open set 
of fibers of some $A \r R_A$ are all strongly F-regular.  
However, this is true when $R$ is Gorenstein.

The point is that for a graded  Gorenstein ring $R$, we may choose a 
homogeneous system  of parameters $\xs$ for $R$ and a socle generator
$z$ for $R/(\xs)R$. We now choose $A$ so that the $\xs$ and the socle generator
$z$ are all defined over $A$. In almost all of the closed fibers $R_{\bar A}$,
the images of the $\xs$ are a system of parameters for the graded 
Gorenstein ring $R_{\bar A}$ and the image of $z$ is a socle generator
for the quotient.  

The fact that all ideals of  $R$  are tightly closed implies that 
$\bar z$ cannot be in  the tight closure $(\bar x_1, \dots, \bar x_d)^*$
in the fiber $R_{\bar A}$, for almost all fibers.
 Because $R$ is Gorenstein, this implies that 
all ideals of $R_{\bar A}$ are tightly closed, so that in particular 
$R_{\bar A}$ is F-split.
 We conclude that  $R_{\bar A}$ is strongly
F-regular on a dense open set of $\maxspec A$. 

Even if we had begun with a $\ZZ$-algebra $A$ over which 
the given system of parameters $\xs$ and the socle generator $z$
 were not defined,
we can still conclude that $R_{\bar A}$ is strongly F-regular on
a non-empty Zariski open set of $\spec A$.  
Simply enlarge $A$ to $A'$ over which they are defined.
Note that $R_{A/{\mu \cap A}} \subset R_{A'/\mu}$ is faithfully flat for
all $\mu \in \maxspec A'$.
 Therefore, because all ideals of $R_{\bar A'} $ are tightly closed (on a 
dense open set), the
same is true of $R_{\bar A}$. In particular, $R_{\bar A}$ is F-split.
Thus  $R$ is D-simple and F-split, and hence strongly F-regular,  
on a dense open set of any finitely generated $\ZZ$-algebra over which $R$ is
defined.
This concludes the proof{\footnote{An alternate way to finish of the proof is to invoke Theorem 2.3.15
of \cite{HH3} to conclude that each of the $R_{\bar A}$ are Gorenstein, 
so that strong F-regularity is equivalent to all ideals being tightly closed
\cite{HH2}. 
We expect that these technical points will eventually be worked out and appear
in the developing manuscript \cite{HH3}.}}.
\end{proof}
\begin{theorems}
Let $S^G$ be the invariant ring for the action of a reductive group
 $G$ on the symmetric algebra $S$
for a finitely dimensional representation  of $G$ of characteristic zero.
Then $R$ is strongly F-regular.
\end{theorems}
\begin{proof}
Because $G$ is reductive, there is a subgroup $H$ of $G$ which is semi-simple
and such that the quotient $G/H$ is an extension of a finite group by a torus.
Note that the quotient group $G/H$ acts on the ring of invariants  $S^H$ for 
the semi-simple group: $\bar g \in G/H$ acts on $f \in S^H$ by
$g\cdot h$ where $g$ is any lifting to $G$ of $\bar g$.
It is easy to verify that $S^G = {(S^H)}^{G/H}$.
Because $H$ is semi-simple, the ring $S^H$ is Gorenstein. Thus by
the preceding lemma, it is strongly F-regular. On the other hand, $G/H$ is 
linearly reductive and thus
the inclusion ${(S^H)}^{G/H} \hookrightarrow S^H$ is split by the Reynolds
operator. This splitting
descends to characteristic $p$ for all $p>0$. Therefore,
because $S^H$ is strongly F-regular in almost all fibers,
 so is its direct summand
$S^G =  {(S^H)}^{G/H}$.
 \end{proof}

The issue of how to keep track of the property of finite F-representation
type in descending to characteristic $p>0$ is wide open. 
\begin{questions}
If $R$ is a  graded direct summand of a polynomial
ring of characteristic zero, does $R$ have Finite F-representation type?
How about if $R$ is a ring of invariants for a semi-simple group acting
linearly on a polynomial ring?
\end{questions}

\ifx\undefined\bysame
\newcommand{\bysame}{\leavevmode\hbox to3em{\hrulefill}\,}
\fi

\end{document}